\documentclass[12pt]{article}

\usepackage[table]{xcolor}
\usepackage{amsmath, amsthm}
\usepackage{amsfonts,amssymb}
\usepackage{hyperref}
\usepackage{authblk}
\usepackage[all]{xy}
\usepackage{enumitem}
\usepackage{stmaryrd}
\usepackage[makeroom]{cancel}
\usepackage{bbm}
\usepackage{fancyvrb}
\usepackage[linesnumbered,lined,boxed]{algorithm2e}
\usepackage{xstring} 
\usepackage{float}
\usepackage{dsfont}
\usepackage{caption}
\usepackage{tikz}
\usetikzlibrary{intersections}

\newif\ifextmath
\extmathfalse

\newcommand{\N}{\mathbb{N}}
\newcommand{\Z}{\mathbb{Z}}
\newcommand{\Q}{\mathbb{Q}}

\newcommand{\C}{\mathbb{C}}
\newcommand{\F}{\mathbb{F}}

\renewcommand{\l}{\mathfrak{l}}

\newcommand{\PP}{\mathbb{P}}
\newcommand{\A}{\mathbb{A}}
\renewcommand{\L}{\mathcal{L}}
\renewcommand{\O}{\mathcal{O}}
\newcommand{\T}{\mathbb{T}}
\newcommand{\G}{\mathbb{G}}

\newcommand{\Gal}{\operatorname{Gal}}
\newcommand{\GQ}{\Gal(\overline \Q / \Q)}
\newcommand{\GL}{\operatorname{GL}}
\newcommand{\SL}{\operatorname{SL}}
\newcommand{\SLZ}{\SL_2(\Z)}
\newcommand{\SLZN}{\SL_2(\Z/N\Z)}
\newcommand{\GLZN}{\GL_2(\Z/N\Z)}
\newcommand{\ZNX}{(\Z/N\Z)^\times}

\newcommand{\Fl}{{\F_\ell}}

\renewcommand{\Im}{\operatorname{Im}}
\newcommand{\Tr}{\operatorname{Tr}}

\newcommand{\disc}{\operatorname{disc}}
\newcommand{\Res}{\operatorname{Res}}
\newcommand{\Aut}{\operatorname{Aut}}
\newcommand{\End}{\operatorname{End}}

\newcommand{\Hom}{\operatorname{Hom}}

\newcommand{\ord}{\operatorname{ord}}

\newcommand{\Frob}{\operatorname{Frob}}

\newcommand{\Pic}{\operatorname{Pic}}

\newcommand{\sgn}{\operatorname{sgn}}

\newcommand{\smat}[4]{\left[ \begin{smallmatrix} #1 & #2 \\ #3 & #4 \end{smallmatrix} \right]}
\newcommand{\mat}[4]{\left[ \begin{matrix} #1 & #2 \\ #3 & #4 \end{matrix} \right]}
\newcommand{\matabcd}{\mat{a}{b}{c}{d}}
\newcommand{\smatabcd}{\smat{a}{b}{c}{d}}
\newcommand{\eps}{\varepsilon}

\newcommand{\mact}{\big\vert}
\newcommand{\Ei}{\mathcal{E}}
\renewcommand{\S}{\mathcal{S}}
\newcommand{\M}{\mathcal{M}}
\newcommand{\newf}{\mathcal{N}}
\newcommand{\W}{\mathcal{W}}

\newcommand{\mfurl}[1]{\StrSubstitute{#1}{.}{/}[\mfslash]}
\newcommand{\mfref}[1]{\mfurl{#1}\href{http://www.lmfdb.org/ModularForm/GL2/Q/holomorphic/\mfslash}{#1}}

\SetKw{True}{True}
\SetKw{False}{False}
\newcommand{\la}{\leftarrow}
\SetKw{FAIL}{FAIL}
\SetKw{Goto}{go to line}

\numberwithin{equation}{subsection}
\numberwithin{figure}{subsection}

\newtheorem{thm}[equation]{Theorem}
\newtheorem{lem}[equation]{Lemma}
\newtheorem{pro}[equation]{Proposition}
\newtheorem{cor}[equation]{Corollary}
\theoremstyle{definition}
\newtheorem{de}[equation]{Definition}
\newtheorem{rk}[equation]{Remark}

\newtheorem{ex}[equation]{Example}

\newcounter{para}[subsection]

\makeatletter
\newcommand{\subjclass}[2][2010]{%
  \let\@oldtitle\@title%
  \gdef\@title{\@oldtitle\footnotetext{#1 \emph{Mathematics subject classification:} #2}}%
}
\newcommand{\keywords}[1]{%
  \let\@@oldtitle\@title%
  \gdef\@title{\@@oldtitle\footnotetext{\emph{Key words and phrases.} #1.}}%
}
\let\c@table\c@equation
\let\c@figure\c@equation
\makeatother
\makeatother

\title{Moduli-friendly Eisenstein series over the~$p$-adics \\ and the computation of \\ modular Galois representations}
\author{Nicolas Mascot\thanks{\href{mailto:mascotn@tcd.ie}{mascotn@tcd.ie}}}
\affil{\scriptsize{Trinity College, Dublin}}

\VerbatimFootnotes

\begin{document}

\maketitle

\begin{abstract}
We show how our~$p$-adic method to compute Galois representations occurring in the torsion of Jacobians of algebraic curves can be adapted to modular curves. The main ingredient is the use of ``moduli-friendly'' Eisenstein series introduced by Makdisi, which allow us to evaluate modular forms at~$p$-adic points of modular curves and dispenses us of the need for equations of modular curves and for~$q$-expansion computations. The resulting algorithm compares very favourably to our complex-analytic method.
\end{abstract}

\vspace{15mm}


\newpage

\section{Introduction}

In this article, given an integer~$k$ and a subgroup~$\Gamma$ of~$\SLZ$ of finite index, we will denote by~$\M_k(\Gamma)$ (resp.~$\S_k(\Gamma)$,~$\Ei_k(\Gamma)$) the space of modular forms (resp. cusp forms, Eisenstein series) of weight~$k$ and level~$\Gamma$.

Let~$f = q + \sum_{n \geqslant 2} a_n q^n \in \S_k\big(\Gamma_1(N)\big)$ be a newform, let~$\l$ be a finite prime of the number field~$\Q(a_n, n \geqslant 2)$, and let~$\rho_{f,\l} : \GQ \longrightarrow \GL_2(\F_\l)$ be the mod~$\l$ Galois representation attached to~$f$. Following Couveignes's and Edixhoven's ideas~\cite{CE11}, we presented in~\cite{algo} and in~\cite{companion} a method to compute~$\rho_{f,\l}$ explicitly, that is to say to find a squarefree polynomial~$F(x) \in \Q[x]$ and a bijection between the roots of~$F(x)$ in~$\overline \Q$ and the nonzero points of the representation space~$\F_\l^2$ such that 
\begin{equation}\text{the Galois action on these roots matches the representation } \rho_{f,\l}. \label{eqn:FGal} \end{equation}

Indeed, these data allow one to determine efficiently the image by~$\rho_{f,\l}$ of Frobenius elements, thanks to the Dokchitsers' method~\cite{Dok}. This method to compute~$\rho_{f,\l}$ is based on complex-analytic geometry and on Makdisi's algorithms~\cite{Mak1},~\cite{Mak2} to compute in Jacobians of curves. It is fairly general, but suffers from some limitations, cf. subsection~\ref{sect:compare} below.

Later on, in~\cite{Hensel}, we presented another method, based on an adaptation of Makdisi's algorithms to a~$p$-adic setting, to compute Galois representations occurring in the torsion of Jacobians of any (not necessarily modular) algebraic curve given by an explicit model (e.g. a plane equation). This method proceeds by computing torsion points over~$\overline \F_p$, and then lifting them~$p$-adically.

The goal of this article is to present an adaptation of this new~$p$-adic method to modular curves, so as to compute representations such as~$\rho_{f,\l}$~$p$-adically. A particularly nice feature of this~$p$-adic approach is that it suppresses the need for plane equations of modular curves, and makes an extremely limited use of~$q$-expansions (cf. section~\ref{sect:Eval}). Besides, it is an occasion to test our~$p$-adic algorithms~\cite{Hensel} in higher genera.

This new approach requires evaluating modular forms at~$p$-adic points of modular curves, which is non-trivial since one cannot use~$q$-expansions for this purpose. We overcome this difficulty by using modular forms introduced in~\cite{MakEis} whose interpretation in terms of the moduli problem parametrised by the modular curve is completely transparent. This is also the reason why we are able to compute in modular Jacobians without requiring equations for the corresponding modular curve. More specifically, the techniques introduced in~\cite{MakEis} make it possible to compute in the Jacobian of a modular curve of level~$N \in \N$ by using only the coordinates of the~$N$-torsion points of a single elliptic curve~$E$ as input data. Moreover, these techniques are completely algebraic, which makes them compatible with base-change and thus usable over~$p$-adic fields (where~$p \nmid 6N$), finite fields (of characteristic coprime to~$6N$), and intermediate objects such as~$\Z/p^e\Z$ with arbitrary~$p$-adic precision~$e \in \N$.

We can thus perform~$p$-adic computations in modular Jacobians with arbitrary finite~$p$-adic accuracy, and hence compute explicitly mod~$\ell$ Galois representations occurring in the~$\ell$-torsion of such Jacobians, and in particular mod~$\ell$ Galois representations attached to eigenforms, thanks to the~$p$-adic method introduced in~\cite{Hensel}.

\bigskip

This article is organised as follows. We begin by recalling the ideas behind our~$p$-adic method~\cite{Hensel} in section~\ref{sect:Hensel}, so as to establish the list of difficulties that we must overcome in order to adapt this method to modular curves. Next, in section~\ref{sect:modcrv}, we gather arithmetic results about modular curves, and in particular about their cusps and the Galois action on them, most of which are probably well known to experts but are unfortunately scattered across the literature. Then, in section~\ref{sect:MakEis}, we recall the definition and some of the properties of Makdisi's moduli-friendly Eisenstein series. After this, we explain in section~\ref{sect:MakMod} how to combine all these ingredients so as to be able to perform~$p$-adic computations in modular Jacobians. The piece of the~$\ell$-torsion of the Jacobian which affords~$\rho_{f,\l}$ can be carved out by using only the action of the Frobenius at~$p$ in most but not all cases, so we show in section~\ref{sect:Tp} how to carve out this piece by using the Hecke operator~$T_p$ instead of the Frobenius. In order to complete the computation of modular Galois representations, it then remains to construct ``evaluation maps'' from the Jacobian to~$\A^1$, which we do in section~\ref{sect:Eval}. Finally, we demonstrate in section~\ref{sect:Examples} that our implementation of the~$p$-adic method presented in this article using~\cite{gp}'s C library outperforms our~\cite{Sage} implementation of the complex-analytic method by a factor ranging from 10 to 100, meaning that computations of moderately ``small'' modular Galois representations now take minutes instead of hours of CPU time, and we explain this difference of performance.

\bigskip

\section{Computing~$p$-adically in Jacobians}\label{sect:Hensel}

\subsection{$p$-adic models of Jacobians}

Let us begin by summarising the kind of data that we need so as to describe a curve in whose Jacobian we want to compute~$p$-adically by using our methods presented in~\cite{Hensel}, which are themselves based on Makdisi's algorithms~\cite{Mak1},~\cite{Mak2}.

%
%
%
%

\begin{de}\label{de:Mak_model}
Let~$C$ be a projective, geometrically non-singular curve of genus~$g$ defined over~$\Q$. Let~$p \in \N$ be a prime number at which~$C$ has good reduction,~$a \in \N$ an integer,~$q=p^a$,~$\Q_q$ the unramified extension of~$\Q_p$ of degree~$a$,~$\Z_q$ its ring of integers,~$\F_q$ its residue field, and finally let~$e \in \N$. A \emph{$p$-adic Makdisi model of~$C$ with residue degree~$a$ and~$p$-adic accuracy~$O(p^e)$} consists of:
\begin{itemize}
\item A choice of a line bundle~$\L$ on~$C$ whose degree~$d_0 = \deg \L$ satisfies
\begin{equation} d_0 \geqslant 2g+1, \label{eqn:d0bound} \end{equation}

\item A choice of points~$P_1, P_2, \cdots, P_{n_Z} \in C(\Z_q/p^e)$, whose number~$n_Z$ satisfies
\begin{equation} n_Z > 5 d_0, \label{eqn:nZbound} \end{equation}
which reduce to pairwise distinct points of~$C(\F_q)$, and which are globally invariant under~$\Frob_p$, as well as the permutation describing the action of~$\Frob_p$ on these points,

\item A choice of a local trivialisation~$t_i : \L \underset{\text{near } P_i}{\simeq} \O_C$ of~$\L$ defined over~$\Q$ (or more generally, over~$\Q_p$) at each of the points~$P_i$, so that we have a Galois-equivariant concept of ``value'' of a global section of~$\L$ at each~$P_i$; we will exclusively use the term ``value'' (with quotation marks) in this sense from now on,

\item A matrix~$V$ of size~$n_Z \times d$ and coefficients in~$\Z_q/p^e$, where
\[ d = \dim H^0(C,\L) = d_0+1-g, \]
whose~$i,j$-entry is the ``value'' in~$\Z_q/p^e$ of~$v_j$ at the point~$P_i$, where the~$v_j$ form a~$\Q_q$-basis of~$H^0(C,\L)$ such that the ``values''~$v_j(P_i)$ lie in~$\Z_q$ and such that the~$v_j$ remain an~$\F_q$-basis of~$H^0(C_{\F_q},\L)$,

\item The local L factor~$L_p(x) \in \Z[x]$ of~$C$ at~$p$, that is to say the numerator of the Zeta function of~$C/\F_p$ reversed so that it is monic and has constant coefficient~$p^g$.
\end{itemize}
\end{de}

\begin{rk}\label{rk:whybounds}
Makdisi's algorithms deal with global sections of powers~$\L^{\otimes n}$ with~$n$ up to~$5$. The bound~\eqref{eqn:nZbound} ensures that such sections are faithfully represented by their ``values'' at the points~$P_i$.

Similarly, the bound~\eqref{eqn:d0bound} ensures that we do avoid complications stemming from dealing with Riemann-Roch spaces attached to divisors of low degree, so that Makdisi's algorithms to compute in Jacobians are valid. For instance, it ensures that the multiplication map
\begin{equation} H^0(C,\L) \otimes H^0(C,\L) \longrightarrow H^0(C,\L^{\otimes 2}) \label{eqn:MakMul} \end{equation}
is surjective, so that we can compute the global sections of~$\L^{\otimes n}$,~$n \leqslant 5$ from the datum of the matrix~$V$.
\end{rk}

\begin{rk}
Note that in particular, such a~$p$-adic Makdisi model of~$C$ does not include an explicit model for~$C$. This is because~\eqref{eqn:d0bound} ensures that~$\L$ is very ample and thus defines a projective embedding of~$C$, whose equations could be read off the kernel of multiplication maps such as~\eqref{eqn:MakMul}. But of course, in order to construct such a~$p$-adic Makdisi model of~$C$, some explicit data about~$C$ must be known, so as to be able to write down the matrix~$V$.
\end{rk}

We show in~\cite{Hensel} that with such a~$p$-adic Makdisi model of~$C$, we can compute the representations of~$\GQ$ afforded in the torsion of the Jacobian~$J$ of~$C$. More precisely, let~$\ell \nmid p$ be a prime, let~$T \subseteq J[\ell]$ be a Galois-submodule, let
\[\rho_T : \GQ \longrightarrow \GL(T)\]
be the mod $\ell$ Galois representation afforded by $T$, and denote by
\[ \chi_T(x) = \det\big(x - \rho_T(\Frob_p)\big) \in \Fl[x] \]
the characteristic polynomial of the Frobenius~$\Frob_p$ acting on~$T$. If
\begin{equation}a \text{ has been chosen so that } \Gal(\overline \Q_q/\Q_q) \text{ acts trivially on } T, \label{eqn:Fq_split_rho}\end{equation}
so that the points of~$T$ are defined over~$\Q_q$, and if
\begin{equation}\chi_T(x) \text{ is coprime mod } \ell \text{ with its cofactor }L_p(x)/\chi_T(x), \label{eqn:chi_mult_1} \end{equation}
so that the datum of~$\chi_T(x)$ determines the submodule~$T \subseteq J[\ell]$ non-ambiguously, then we can compute~$\rho_T$ as follows:

\begin{figure}[H]
\fbox{\parbox{\textwidth}{
\begin{enumerate}
\item Determine the order~$N = \# J(\F_q)$ as~$N = \Res(L_p(x),x^a-1)$, and factor it as~$N = \ell^v M$ where~$\ell \nmid M$,
\item \label{algo:Strategy_Hensel_randtors} Take a random point~$x \in J(\F_q)$, multiply it by~$M$ so as to get an~$\ell$-power torsion point, then repeatedly by~$\ell$ so as to get an~$\ell$-torsion point, and finally apply~$\frac{L_p}{\chi}(\Frob_p)$ to it so as to project it onto~$T$,
\item Repeat this process so as to obtain an~$\Fl$-basis of~$T$ over~$\F_q$,
\item \label{algo:Strategy_Hensel:Lift} Lift this basis of~$T$ to~$J(\Z_q/p^e)[\ell]$,
\item Use Makdisi's algorithms over $\Z_q/p^e$ to compute all the points of~$T$ over $\Z_q/p^e$, by forming all linear combinations from this basis,
\item Construct a rational map~$\alpha : J \dashrightarrow \A^1$ defined over~$\Q$,
\item Compute~$F(x) = \prod_{t \in T} \big(x - \alpha(t) \big) \in \Z_q/p^e[x]$, and identify it as an element of~$\Q[x]$.
\end{enumerate}
}}
\caption{$p$-adic computation of Galois representations found in Jacobians}
\label{algo:Strategy_Hensel}
\end{figure}

Indeed, if~$\alpha$ is sufficiently generic to be injective on~$T$, then the values~$\alpha(t)$ are permuted by~$\GQ$ in a way matching the points of~$T$, so the polynomial~$F(x)$ satisfies~\eqref{eqn:FGal}.

\begin{rk}
We actually get a mod~$p^e$ approximation of~$F(x)$, so we need the accuracy parameter~$e$ to be large enough so as to  be able to identify~$F(x) \in \Q[x]$. Naturally, higher values of~$e$ will also work, but slow the computation down. As in~\cite{Hensel}, we do not have a clear recipe for the ideal value of~$e$, and we proceed mostly by trial-and-error. Therefore, the results which we get are not rigorously proved to be correct; however, in the case of Galois representations attached to modular forms, they can be rigorously certified using the methods presented in~\cite{certif}. In what follows, we will not concern ourselves with this aspect anymore, and simply assume that the value of~$e$ has been set somehow.
\end{rk}


In order to adapt this method to compute Galois representations to modular curves, we must construct a~$p$-adic Makdisi model of these modular curves. It is natural to choose the line bundle~$\L$ so that its sections are modular forms, but then we need to be able to evaluate these modular forms at~$p$-adic points of the modular curve so as to be able to write down the matrix~$V$. We explain how this can be done efficiently in the rest of the article, but before that, we introduce two improvements to the method~\cite{Hensel} which we should have included in~\cite{Hensel} and will be useful for our purpose later on.

\begin{rk}
Strategy~\ref{algo:Strategy_Hensel} assumes that we can find a good prime~$p$ satisfying~\eqref{eqn:chi_mult_1}. This is very often possible, but not always, as demonstrated by example~\ref{ex:Frobp_cannot_cut} below. We explore a remedy to this unpleasant situation in section~\ref{sect:Tp}.
\end{rk}

\subsection{Automorphisms and Frobenius}\label{subs:JAuts}

By Riemann-Roch and assumption~\eqref{eqn:d0bound}, every point~$ x \in J = \Pic^0(C)$ is represented by the line bundle~$\L(-D_x)$ for some (non unique) effective divisor~$D_x$ on~$C$ of degree~$d_0 = \deg \L$. A point~$x \in J(\Q_q)$ may thus be represented by the matrix
\[ W_{D_x} = \Big( w_j(P_i) \Big)_{\begin{array}{ll}\scriptstyle i \leqslant n_Z \\ \scriptstyle j \leqslant d_0+1-g \end{array}} \]
where the~$P_i$ are as above and the~$w_j$ form a~$\Q_q$-basis of the space of global sections of~$\L^{\otimes 2}(-D_x)$ chosen so that the ``values''~$w_j(P_i)$ lie in~$\Z_q$ and that the~$w_j$ still determine a basis of this section space over~$\F_q$. As explained throughout section~2 of~\cite{Hensel}, this mode of representation of points of~$J$ allow us to perform~$p$-adic computations in~$J$ with a~$p$-adic Makdisi model of~$C$ thanks to Makdisi's algorithms; but naturally, as explicit computations (and in particular our~$p$-adic Makdisi model) can only involve a finite~$p$-adic accuracy, our algorithms only deal with points of~$J(\Z_q/p^e)$, which are internally represented by a matrix~$W_{D_x}$ defined as above but whose entries lie in~$\Z_q/p^e$.

Since our local trivialisations of~$\L$ are defined over~$\Q_p$, we have
\[ w_j^{\Frob_p}(P_i) = \left(w_j\big(P_i^{\Frob_p^{-1}}\big)\right)^{\Frob_p} \]
for all~$i$ and~$j$. As explained in~\cite[2.2.5]{Hensel}, this means that given a matrix~$W_{D_x}$ representing a point~$x \in J(\Z_q/p^e)$ as above, we may obtain the matrix~$W_{D_x^{\Frob_p}}$ representing the point~$x^{\Frob_p}$ of~$J(\Z_q/p^e)$ by applying~$\Frob_p$ to the entries of~$W_{D_x}$ and permuting its rows by the inverse of the permutation induced by~$\Frob_p$ on the points~$P_i \in C(\Z_q/p^e)$, which we can do since this permutation is recorded as part of the~$p$-adic Makdisi model. Note that compared to the group law in~$J(\Z_q/p^e)$, which involves linear algebra on matrices~$W_D$, this process is almost instantaneous.

\bigskip

Suppose now that we have an automorphism~$\varphi \in \Aut(C)$ of~$C$ which is defined over~$\Q$ (or more generally, over~$\Q_p$). As explained in~\cite[6.2]{DS}, it extends by linearity to a map~$\varphi_*$ on divisors of~$C$, which in turn induces an automorphism of~$J$ which we also denote by~$\varphi_*$, because the norm map~$N_\varphi : \O_C \longrightarrow \O_C$, which takes a section~$f$ to~$f \circ \varphi^{-1}$, satisfies~$\varphi_*\big((f)\big) = \big(N_{\varphi}(f)\big)$. Suppose that during the construction of the~$p$-adic Makdisi model of~$C$, we have chosen a line bundle~$\L$ which satisfies~$\varphi_* \L = \L$, points~$P_i \in C(\Q_q)$ which are globally invariant under~$\varphi$, that we have recorded the permutation~$\sigma_\varphi$ defined by~$\varphi(P_i) = P_{\sigma_\varphi(i)}$, and that our local trivialisations~$t_i$ are compatible with~$\varphi$, in that
\[ \xymatrix{ \L \ar[r]^{t_i} \ar[d]_{\varphi_*} & \O_C \ar[d]^{\varphi_*} \\
\L \ar[r]_{t_{\sigma(i)}} & \O_C} \]
commutes for all~$i$. Then, by the same line of ideas as for~$\Frob_p$, we may instantaneously apply~$\varphi_*$ to a matrix~$W_{D_x}$ representing a point~$x \in J$: all we have to do is permute its rows by~$\sigma_\phi^{-1}$.

Indeed, on the one hand, as~$x \in J$ is represented by~$\L(-D_x)$, its image~$\varphi_*(x)$ is represented by~$\L(-\varphi_*(D_x))$ since~$\L$ is invariant by~$\varphi$; and on the other hand, each column of~$W_{D_x}$ consists of the vector of ``values''~$s(P_i)$ of a global section~$s$ of~$\L^{\otimes 2}(-D_x)$, so permuting its entries by~$\sigma_{\varphi}^{-1}$ yields the vector of ``values'' of a section~$s'$ such that~$s'(P_{\sigma(i)}) = s(P_i)$, so that indeed~$s' = N_\varphi(s)$.

We will use this idea later on in this article with~$C$ a modular curve~$X_H(N)$ and~$\varphi$ a diamond operator~$\langle d \rangle$ (cf. subsection~\ref{subsect:modcrv} below for definitions).

\subsection{Fast exponentiation using cyclotomic \\ polynomials and the Frobenius}\label{sect:cycloexp}

In the notation of strategy~\ref{algo:Strategy_Hensel}, we typically have~$M \approx \# J(\F_q) \approx q^g$. Therefore, for large genus~$g$, on the one hand~$M$ is quite large, especially as~\eqref{eqn:Fq_split_rho} usually imposes that~$q$ is in the thousands or even millions; and on the other hand, performing one addition in~$J$ using Makdisi's algorithms relies on linear algebra of size~$O(g)$ which is rather costly. Thus, even though we use fast exponentiation, multiplication by~$M$ in~$J(\F_q)$, which is a required step to generate~$\ell$-torsion points, can take a significant amount of time. However, we have seen in the previous subsection that applying the Frobenius~$\Frob_p$ is almost instantaneous, so it is natural to try to use the action of Frobenius in order to speed up this multiplication-by-$M$ step. For this, we begin by establishing the following result:

\begin{lem}\label{lem:decomp}
Let~$G$ be a finite Abelian group, and let~$\phi : G \rightarrow G$ be an endomorphism. View~$G$ as a~$\Z[x]$-module with~$x$ acting as~$\phi$. Suppose we know a monic polynomial~$F(x) \in \Z[x]$ such that~$F(\phi)=0$, and that~$F$ factors in~$\Z[x]$ as~$F=AB$ with~$\Res(A,B)$ coprime to~$\# G$ (in particular,~$A$ and~$B$ must be coprime in~$\Q[x]$). Then~$G$ decomposes as~$G[A] \times G[B]$.
\end{lem}

\begin{proof}
%
Let $U_1,V_1 \in \Z[x]$ be such that $U_1A+V_1B = \Res(A,B)$. Since $\Res(A,B)$ and $\# G$ are coprime, we can find $m \in \Z$ such that $m \Res(A,B) \equiv 1 \bmod \#G$; multiplying $U_1$ and $V_1$ by $m$ thus yields $U,V \in \Z[x]$ such that $UA+VB=1$.
It is then clear that the maps
\[ \begin{array}{ccc}
G & \longleftrightarrow & G[A] \times G[B] \\
g & \longmapsto & (VBg, \ UAg) \\
a+b & \longmapsfrom & (a, \ b) \\
\end{array} \]
are inverses of each other.
\end{proof}

Suppose now that~$\ell \nmid a$; if one has chosen~$a$ minimal for the points of~$T$ to be defined over~$\F_q = \F_{p^a}$, which is what we will do in practice, then this is equivalent to saying that~$\rho_T(\Frob_p)$, which has order $a$, is semisimple. Take~$G$ to be the part of~$J(\F_q)$ coprime to~$a$ (which has thus the same~$\ell$-torsion as~$J(\F_q)$),~$\phi = \Frob_p : G \rightarrow G$, and~$F(x)=x^a-1$, which factors over~$\Z$ into the cyclotomic polynomials
\[ x^a-1 = \prod_{d \mid a} \Phi_d(x). \]
Since~$\# G$ is prime to~$\disc(x^a-1) = \pm a^a$ by construction, an iterated use of lemma~\ref{lem:decomp} yields the decomposition
\[ G = \bigoplus_{d \mid n} G[\Phi_d]. \]
The point is that
\[ \# G[\Phi_d] \approx N_d \overset{\text{def}}{=} J[\Phi_d] = \Res(L_p,\Phi_d) = \prod_{L_p(\alpha)=0} \Phi_d(\alpha) \approx q^{g \varphi(d)} \]
so each of these factors is typically much smaller than~$J(\F_q)$.

\bigskip

We can thus obtain~$\ell$-torsion points as follows:

\begin{figure}[H]
\fbox{\parbox{\textwidth}{
\begin{enumerate}
\item Pick~$d \mid a$ such that~$\ell \mid N_d$ (there will be at least one).
\item Write~$N_d= \ell^{v_d} M_d$ with~$M_d$ coprime to~$N$.
\item Take a random~$x \in J(\F_{q})$. Apply~$\frac{x^a-1}{\Phi_d(x)}(\Frob_p)$ to it, multiply the result by~$M_d$, and then repeatedly by~$\ell$ until we get 0. Return the last nonzero point.
\end{enumerate}
}}
\caption{Generating torsion points thanks to cyclotomic polynomials}
\label{algo:Strategy_CycloExp}
\end{figure}

This method is valid since as~$\ell \nmid a$ by assumption,~$M_d$ is divisible by~$[J(\F_{q}):G]$, so that multiplication by~$M_d$ includes the effect of projecting from~$J(\F_{q})$ to~$G$. Its advantage is that the number of required operations in~$J(\F_q)$ is approximately \[ \lg N_d \approx g \varphi(d) \lg p, \] compared with \[\lg N \approx g a \lg p\] with strategy~\eqref{algo:Strategy_Hensel}. Indeed, typically the cofactor~$\frac{x^a-1}{\Phi_d(x)}$ has few nonzero coefficients, and these coefficients are usually~$\pm1$, so applying~$\frac{x^a-1}{\Phi_d(x)}(\Frob_p)$ requires few operations in~$J(\F_q)$ and thus takes negligible time since applying~$\Frob_p$ is instantaneous.

\begin{ex}
Let~$C$ be the modular curve~$X_1(13)$, which has genus 2, and let~$J$ be its Jacobian. Suppose we want to generate~$\ell$-torsion points of~$J$ where~$\ell=29$ for example. Take~$p=191$; using formula~\eqref{eqn:Lp_mod_res} below and~\cite[proposition 5.1]{Hensel}, we find that the smallest~$a \in \N$ such that~$J[\ell]$ is defined over~$\F_{p^a}$ is~$a=12$ (which is why we chose this~$p$, as other values of~$p$ typically require~$a$ to be in the hundreds if not more).

We have
\[ \# J(\F_{p^a}) = \ell^4 M \]
where
\[ \lg M \approx 162, \]
so if we use the method presented in strategy~\eqref{algo:Strategy_Hensel}, then we need to perform about~$200$ additions in~$J(\F_{p^a})$ in order to obtain an~$\ell$-torsion point.

In comparison, if we take~$d=12$, we find
\[ N_{12} = \ell^2 M_{12} \]
where
\[ \lg M_{12} \approx 60, \]
so we can produce an~$\ell$-torsion point with less than~$100$ additions in~$J(\F_{p^a})$ by using strategy~\eqref{algo:Strategy_CycloExp}, even taking into account the operations required to multiply by~$\frac{x^{n}-1}{\Phi_{d}(x)} = x^8 + x^6 - x^2 - 1$.

Similarly, for~$d=3$ we have
\[ N_3 = \ell^2 M_3 \]
so this is the only~$d$ the rest of~$J[\ell]$ comes from, and we have
\[ \lg M_{3} \approx 20 \]
only. However, this produces~$\ell$-torsion points defined over~$\F_{p^3}$, so if we want to get all of~$J[\ell]$, then we need to generate points using~$d=12$ as well.
\end{ex}

\begin{rk}
In our case, we do not only want points of~$J[\ell]$, but actually points in the piece~$T_\chi$ of~$J[\ell]$ where~$\Frob_p$ acts with characteristic polynomial~$\chi(x)$. This means that in strategy~\eqref{algo:Strategy_CycloExp}, we should only consider the~$d \mid a$ such that~$\Phi_d(x)$ has a nontrivial common factor mod~$\ell$ with~$\chi(x)$.
\end{rk}


\section{Reminders on modular curves and their cusps}\label{sect:modcrv}

\subsection{Classical congruence subgroups and their moduli problems}\label{subsect:modcrv}

Let~$N \in \N$, and define as usual
\[ \Gamma(N) = \left\{ \gamma \in \SL_2(\Z) \ \vert \ \gamma \equiv 1 \bmod N \right\}, \]
\[ \Gamma_0(N) = \left\{ \gamma = \matabcd \in \SL_2(\Z) \ \vert \ c \equiv 0 \bmod N \right\}, \]
\[ \Gamma_1(N) = \left\{ \gamma = \matabcd \in \Gamma_0(N) \ \vert \ a \equiv d \equiv 1 \bmod N \right\}, \]
and more generally, given a subgroup~$H \leqslant (\Z/N\Z)^\times$,
\[ \Gamma_H(N) = \left\{ \gamma = \matabcd \in \Gamma_0(N) \ \vert \ a,d \bmod N \in H \right\}. \]
Denote the corresponding modular curves by~$X(N)$,~$X_0(N)$,~$X_1(N)$, and~$X_H(N)$. Note that since~$-1 \in \SL_2(\Z)$ acts trivially on the upper-half plane, we have~$X_H(N) = X_{\langle H,-1\rangle}(N)$ where~$\langle H,-1\rangle$ denotes the subgroup of~$(\Z/N\Z)^\times$ generated by~$H$ and~$-1$, so we will restrict our attention to the subgroups~$H \leqslant (\Z/N\Z)^\times$ that contain~$-1$. 

\bigskip

Let us briefly recall the interest of these modular curves, and use this occasion to fix some notation and conventions which we will use throughout the rest of this article. Informally speaking, the curve~$X(N)$ parametrises the pairs~$(E,\beta)$ up to isomorphism, where~$E$ is an elliptic curve and
 \[ \beta: (\Z/N\Z)^2  \overset{\sim}{\longrightarrow} E[N] \]
is an isomorphism mapping the standard basis~$[1,0], [0,1]$ of~$(\Z/N\Z)^2$ to points~$P,Q \in E[N]$ such that the Weil paring~$e_N(P,Q)$ is a fixed primitive~$N$-th root of 1.

We pause here to mention that we normalise the Weil-pairing as in~\cite[7.4]{DS}, so that~\cite[1.3]{DS} for any~$\omega_1, \omega_2 \in \C$ such that~$\Im \frac{\omega_1}{\omega_2} > 0$, we have
\begin{equation} e_N(\omega_1 / N,\omega_2 / N) = e^{2\pi i /N} \label{eqn:normalisation_Weil} \end{equation}
on the elliptic curve~$\C/(\Z \omega_1 \oplus\Z\omega_2)$; beware that some authors (and~\cite{gp}) use the opposite normalisation, namely~$e_N(\omega_2 / N,\omega_1 / N) = e^{2\pi i /N}$. This choice of normalisation will matter later (cf. theorem~\ref{thm:qexp} below). 

We will always view the elements of~$(\Z/N\Z)^2$ as \emph{row} vectors; we then have a left action of~$\SLZN$ on~$X(N)$ defined by 
\[ \gamma \cdot (E,\beta) = (E,\beta_\gamma), \]
where~$\gamma \in \SLZN$ and~$\beta_\gamma$ is the isomorphism between~$(\Z/N\Z)^2$ and~$E[N]$ taking~$v \in (\Z/N\Z)^2$ to~$\beta(v \gamma)$. It follows that each fibre of the projection map~$X(N) \rightarrow X(1)$ at an elliptic curve~$E$ having no automorphisms other than~$\pm1$ is a torsor under~$\SLZN / \pm1$.

Similarly, the curve~$X_1(N)$ parametrises isomorphism classes of pairs~$(E,Q)$, where~$E$ is an elliptic curve and~$Q \in E$ is a point of exact order~$N$; more generally,~$X_H(N)$ paramatrises isomorphism classes of pairs~$(E,H\cdot Q)$ of elliptic curves equipped with a point of order~$N$ up to multiplication by~$H$. The projection map from~$X(N)$ to~$X_H(N)$ is given by 
\begin{equation} \begin{array}{ccl} X(N) & \longrightarrow & X_H(N) \\ (E,\beta) &  \longmapsto & \big(E,H \cdot \beta([0,1]) \big), \end{array} \label{eqn:proj_XN_XH} \end{equation}
so that given~$\gamma, \gamma' \in \SLZN$, the points~$(E,\beta_\gamma)$ and~$(E,\beta_{\gamma'})$ of~$X(N)$ project to the same point of~$X_H(N)$ iff.~$\gamma$ and~$\gamma'$ have the same \emph{bottom row} up to scaling by~$H$ (remember that we are assuming that~$-1 \in H$). It follows that given an elliptic curve~$E$ such that~$\Aut(E)$ is reduced to~$\{\pm1\}$ and an isomorphism~$\beta:(\Z/N\Z)^2 \simeq E[N]$, we have a bijection
\begin{equation} \begin{array}{ccc} \text{Primitive vectors of }(\Z/N\Z)^2 \text{ up to }H & \longleftrightarrow & \text{Fibre of } X_H(N) \rightarrow X(1) \text{ at } E \\ H \cdot [c,d] & \longmapsto & \big(E,H \cdot \beta_\gamma([0,1])\big) = \big(E,H \cdot \beta([c,d])\big) \end{array} \label{eqn:fibre_XH} \end{equation}
where~$\gamma$ denotes any element of~$\SLZN$ whose bottom row is~$[c,d]$.

\bigskip

The group~$\Gamma_1(N)$ is normal in~$\Gamma_H(N)$, and we have the isomorphism
\begin{equation} \begin{array}{ccc} \Gamma_H(N) / \Gamma_1(N) & \overset{\sim}{\longrightarrow} & \ZNX / H \\ \matabcd & \longmapsto & d. \end{array} \label{eqn:GammaHmod1} \end{equation}
Given~$y \in \ZNX / H$, we may thus define the \emph{diamond operator}~$\langle y \rangle$ as the automorphism of~$X_H(N)$ which acts as the inverse image of~$y$ by~\eqref{eqn:GammaHmod1}; in other words, under the moduli point of view, it takes the pair~$(E,H \cdot Q)$ to the pair~$(E,y \, H \cdot Q)$, and under~\eqref{eqn:fibre_XH}, it corresponds to
\begin{equation} H \cdot [c,d] \longmapsto H \cdot [yc,yd]. \label{eqn:Diam_on_fibre} \end{equation}

\bigskip

Let~$\mu_N \subset \overline \Q$ denote the group of~$N$-th roots of 1, and identify the Galois group~$\Gal(\Q(\mu_N)/\Q)$ of the~$N$-th cyclotomic field with~$\ZNX$ via
\begin{equation} \begin{array}{ccl} \ZNX & \overset{\sim}{\longrightarrow} & \Gal(\Q(\mu_N)/\Q) \\ x & \longmapsto & \sigma_x : (\zeta \mapsto \zeta^x) \text{ for all } \zeta \in \mu_N. \end{array} \label{eqn:Gal_cyclo} \end{equation}
The moduli interpretation of~$X(N)$ (resp.\ of~$X_1(N)$,~$X_0(N)$, and more generally~$X_H(N)$) makes sense over~$\Q(\mu_N)$ (resp.\ over~$\Q$), so this curve admits a model over~$\Q(\mu_N)$ (resp.\ over~$\Q$) which is compatible with this moduli interpretation, so that in particular the diamond operators are defined over~$\Q$. For what follows, we must describe precisely such a model.

As in~\cite[6.2]{Shimura}, consider the function field
\[ F_N = \Q(j, f_0^v \ \vert \ 0 \neq v  \in (\Z/N\Z)^2), \]
where for each nonzero~$v = (c_v, d_v) \in (\Z/N\Z)^2$, the modular function~$f_0^v$ is defined on the upper-half plane by
\[ f_0^v(\tau) = \frac{G_4(\tau)}{G_6(\tau)}\wp_\tau\left(\frac{c_v \tau+d_v}N\right) \]
where
\[ G_k(\tau) = \sum_{\substack{m,n \in \Z \\ (m,n) \neq (0,0)}} \frac1{(m\tau+n)^k} \]
and~$\wp_\tau$ is the Weierstrass~$\wp$ function attached to the lattice spanned by~$\tau$ and 1. Then~\cite[6.2]{Shimura} we have
\[ F_N \cap \overline \Q = \Q(\mu_N), \]
so~$F_N$ provides us with a model of~$X(N)$ over~$\Q(\mu_N)$.

As $\wp_\tau(z)$ is an even function of $z$ for all $\tau$, we have $f_0^{-v} = f_0^v$ for all $v$, whence a natural right action of~$G_N = \GL_2(\Z/N\Z) / \pm 1$ on~$F_N$ defined by
\[ f_0^v \cdot \gamma = f_0^{v\gamma} \quad (v \in (\Z/N\Z)^2, \gamma \in G_N), \]
making~$F_N$ a Galois extension of~$F_N^{G_N} = \Q(j)$ with Galois group~$G_N$. Furthermore, each~$\gamma \in G_N$ restricts to~$\sigma_{\det \gamma} \in \Gal(\Q(\mu_N)/\Q)$ on~$\Q(\mu_N) = F_N \cap \overline \Q$. Each subgroup~$U \leqslant G_N$ hence corresponds to the function field~$F_N^U$ of a quotient of~$X(N)$ defined over the subfield~$\Q(\mu_N)^{\det(U)}$ of~$\Q(\mu_N)$. In view of~\eqref{eqn:proj_XN_XH}, given a subgroup~$H \leqslant (\Z/N\Z)^\times$ containing~$-1$, we may thus set~\cite[7.7]{DS}
\begin{equation} \Q\big(X_H(N)\big) = F_N^{U_H}, \quad \text{where} \quad  U_H =  \smat{*}{*}{0}{H} = \left\{ \pm \smatabcd \in G_N \ \vert \ c=0, d \in H \right\}, \label{eqn:fnfield_XH} \end{equation}
thus fixing a model for~$X_H(N)$ over~$\Q$ for each such~$H$, and in particular for~$X_0(N)$ and~$X_1(N)$. In particular, this means that for~$N \geqslant 5$,~$X_1(N)$ is the moduli space for elliptic curves~$E$ equipped with a torsion point of exact order~$N$, or equivalently with an embedding~$\Z/N\Z \hookrightarrow E[N]$, as opposed to an embedding~$\mu_N \hookrightarrow E[N]$; cf.~\cite[9.3]{DI} and Example~\ref{ex:cusps_X1} below.

\subsection{Modular Galois representations in modular Jacobians}\label{sect:rho_in_XH}

Let us define a \emph{newform} of level~$\Gamma_H(N)$ as a newform of level~$\Gamma_1(N)$ whose nebentypus ~$\eps \colon (\Z/N\Z)^\times \to \overline \Q^\times$ satisfies~$H \leqslant \ker \eps$. Let~$\newf_k(\Gamma_H(N))$ denote the finite set of newforms of weight~$k$ and level~$\Gamma_H(N)$. This set is acted on by~$\GQ$ via the coefficients of~$q$-expansions at the cusp~$\infty$: if~$f(q)=\sum_{n=1}^\infty a_n(f) q^n$ has nebentypus~$\eps$, then the coefficients~$\{a_n(f)\}_n$ are algebraic integers, and for~$\tau \in \GQ$ we have~$\tau f \in \newf_k(\Gamma_H(N))$ with~$q$-expansion~$(\tau f)(q) = \sum_{n} \tau(a_n(f)) q^n$ and nebentypus~$\tau \circ \eps$. Denote the Galois orbit of~$f$ by~$[f]$, and let
\[ \GQ \backslash \newf_k(\Gamma_H(N)) \]
be the set of such Galois orbits. Then the Jacobian~$J_H(N)$ of~$X_H(N)$ decomposes up to isogeny over~$\Q$ as
\begin{equation}
J_H(N) \sim \prod_{M \mid N} \prod_{[f] \in \GQ \backslash \newf_2( \Gamma_{H_M}(M))} A_{[f]}^{\sigma_0(N/M)}, \label{eqn:decomp_mod_jac}
\end{equation}
where~$\sigma_0(n) = \sum_{d \mid n} 1$ is the number of divisors of~$n$,~$H_M \leqslant (\Z/M\Z)^\times$ denotes the image of~$H$ in~$(\Z/M\Z)^\times$, and for each~$[f]$, the Abelian variety~$A_{[f]}$ is simple over~$\Q$ of dimension~$\dim A_{[f]} = [\Q\big(a_n(f) \ \vert \ n \geqslant 2\big) : \Q]$. Roughly speaking,~$A_{[f]}$ can be thought of as the piece of~$J_H(N)$ where the Hecke algebra~$\T$ of weight~$2$ and level~$\Gamma_H(N)$ acts with the eigenvalue system of~$f$; more precisely,~$A_{[f]}$ is defined by
\[ A_{[f]} = J_H(N) / I_{[f]} J_H(N), \]
where~$I_{[f]} = \{ T \in \T \ \vert \ Tf = 0\}$ is the annihilator of~$f$ under the Hecke algebra.

\bigskip

Let~$f$ be a newform of weight~$k$, level~$N$, and nebentypus~$\eps_f$; let~$K_f$ be the number field~$\Q(a_n(f) \ \vert \ n \geqslant 2)$, which contains the values of~$\eps_f$, let~$\l$ be a finite prime of~$K_f$, and finally let~$\ell \in \N$ be the prime below~$\l$. Suppose we wish to compute the mod~$\l$ Galois representation~$\rho_{f,\l}$ attached to~$f$. By~\cite[2.1]{Ribet}, this representation is also attached to a form whose is level coprime to~$\ell$, so we assume that~$\ell \nmid N$ from now on. Similarly, by theorem 2.7 of~\cite{RS}, up to twist by the mod~$\ell$ cyclotomic character, this representation is also attached to a form of the same level and of weight comprised between~$2$ and~$\ell+1$, so we suppose~$2 \leqslant k \leqslant \ell+1$ from now on.

If~$k=2$, then~$\rho_{f,\l}$ occurs in the~$\ell$-torsion of the Jacobian~$J_1(N)$ of~$X_1(N)$. Else, recall~\cite[p.178]{RS} that there exists an eigenform~$f_2$ of weight 2 but level~$\ell N$ and a prime~$\l_2$ of~$K_{f_2}$ above~$\ell$ such that
\begin{equation} \rho_{f_2,\l_2} \sim \rho_{f,\l}. \label{eqn:rho_equiv_wt2} \end{equation}
We thus set
\begin{equation} N' = \left\{ \begin{array}{ll} N, & k=2 \\ \ell N, & k >2,\end{array} \right. \label{eqn:def_N'} \end{equation}
so that~$\rho_{f,\l}$ occurs in~$J_1(N')[\ell]$; also define~$f_2 = f$ if~$k=2$. Then~$\rho_{f,\l}$ actually occurs in~$A_{[f_2]}[\ell]$, so in view of~\eqref{eqn:decomp_mod_jac},~$\rho_{f,\l}$ occurs in~$J_H(N')$ provided that~$H \leqslant \ker \eps_2$, where~$\eps_2$ is the nebentypus of~$f_2$. Taking~$H = \ker \eps_2$, we thus get a modular curve whose Jacobian contains~$\rho_{f,\l}$, but whose genus is (hopefully) smaller than that of~$X_1(N')$, making explicit computations with it more efficient. This is our reason for introducing the modular curves~$X_H(N)$; this idea originates from~\cite[4.1]{GammaH}. More explicitly,~\eqref{eqn:rho_equiv_wt2} implies that
\[ x^{k-1} \eps_f(x) \bmod \l = x \eps_2(x) \bmod \l_2 \text{ for all } x \in (\Z/N'\Z)^\times, \]
so that we take
\begin{equation} H = \ker \eps_2 = \{ x \in (\Z/N'\Z)^\times \ \vert \ x^{k-2} \eps_f(x) = 1 \bmod \l \}. \label{eqn:def_H} \end{equation}

\begin{rk}
Naturally, in many cases,~$H$ is a very small subgroup of~$(\Z/N'\Z)^\times$, so that the genus of~$X_H(N')$ is the same as, or not much smaller than, that of~$X_1(N')$. However, there are also cases when the genus of~$X_H(N')$ is dramatically smaller than that of~$X_1(N')$, which makes it possible to compute Galois representations that would otherwise be out of reach, cf.~\cite{companion} for some examples.
\end{rk}

\begin{rk}
In principle, it would be even better to compute~$\rho_{f,\l}$ directly in the Abelian variety~$A_{[f_2]}$, but the author only knows how to compute with Jacobians.
\end{rk}

In order to construct a~$p$-adic Makdisi model for~$X_H(N')$, we will in particular need to determine the local L factor of~$X_H(N')$ at a prime~$p \nmid N'$. For this, we suppose that we can compute the set of Galois orbits of mod~$N'$ Dirichlet characters
\[ \chi : (\Z/N'\Z)^\times \rightarrow \Q[t]/\Phi_{\ord \chi}(t), \]
where~$\Phi_n(t) \in \Z[t]$ denotes the~$n$-th cyclotomic polynomial, and that for each such orbit, we can compute the matrix of the Hecke operator~$T_p$ with respect to some~$\Q[t]/\Phi_{\ord \chi}(t)$-basis of the space of cusp forms of level~$N'$, weight~$2$, and nebentypus~$\chi$; for instance, this is possible using~\cite{gp}. Then, in view of the decomposition~\eqref{eqn:decomp_mod_jac}, we have
\begin{equation} L_p\big(X_H(N')\big) = \hspace{-1cm} \prod_{\substack{\chi \bmod \GQ \\ \chi : (\Z/N'\Z)^\times \rightarrow \Q[t]/\Phi_{\ord \chi}(t) \\ \ker \chi \geqslant H}} \hspace{-1cm} \Res_t\big( \Phi_{\ord \chi}(t),\Res_y(x^2-yx+p\chi(p), \det(y1-T_p \vert_{\S_2(N',\chi)}) \big). \label{eqn:Lp_mod_res} \end{equation}
In particular, we also recover the genus of~$X_H(N')$ as half the degree of this polynomial.

\subsection{The cusps of~$X_H(N)$}

\subsubsection*{Moduli interpretation and Galois action}

Let~$N \in \N$. Recall that the \emph{N\'eron~$N$-gon} is the variety~$C_N$ obtained by gluing~$N$ copies of~$\PP^1$ indexed by~$\Z/N\Z$ by the relation 
\[ (\infty,i) \sim (0, i+1) \]
for all~$i \in \Z/N\Z$. Its regular locus is thus the group variety
\[ C_N^{\text{reg}} = \G_m \times \Z/N\Z. \]
A \emph{morphism} between such~$N$-gons is an algebraic variety morphism inducing a group variety morphism on the regular locus.

Whereas the non-cuspidal points of the modular curve~$X_1(N)$ correspond to isomorphism classes of pairs formed by an elliptic curve~$E$ and a torsion point of~$E$ of exact order~$N$, the cusps of~$X_1(N)$ correspond to isomorphism classes of pairs formed by a N\'eron~$n$-gons (for some~$n \in \N$) equipped with a torsion point of~$C_n^{\text{reg}}$ of exact order~$N$ whose multiples meet every component, cf.~\cite[9.3]{DI}. Such a pair is thus of the form
\[ \big( C_n, (\zeta,i) \big) \]
where~$n \mid N$,~$\zeta \in \mu_N$, and~$i \in (\Z/n\Z)^\times$ are such that the gcd of~$N$, the order of~$\zeta$, and the order of~$i$ is~$1$; in particular, it its defined over~$\Q(\zeta) \subseteq \Q(\mu_N)$. In order to understand the cusps of~$X_1(N)$, and in particular how thy are permuted by~$\GQ$, we must therefore classify such pairs up to isomorphism, and in particular determine the automorphisms of~$C_n$.

First of all, we have canonically
\[ \End(C_n^{\text{reg}}) = \End(\G_m \times \Z/n\Z) = \mat{\End(\G_m)}{\Hom(\Z/n\Z,\G_m)}{\Hom(\G_m, \Z/n\Z)}{\End(\Z/n\Z)} =  \mat{\Z}{\mu_n}{0}{\Z/n\Z} \]
acting by 
\[ \mat{m}{\zeta}{0}{j} \left[ \begin{matrix} x \\ i \end{matrix} \right]= \left[ \begin{matrix} \zeta^i x^m \\ ji \end{matrix} \right] \quad (x \in \G_m, i \in \Z/n\Z).  \]
Therefore,
\[ \Aut(C_n^{\text{reg}}) =  \mat{\pm1}{\mu_n}{0}{( \Z/n\Z )^\times}. \]
Finally, an automorphism of~$C_n^{\text{reg}}$ extends to an automorphism of~$C_n$ iff. it respects the gluing condition~$(\infty,i) \sim (0, i+1)$, which translates into the condition~$j=m$. Therefore,
\begin{equation} \Aut(C_n) = \pm \mat{1}{\mu_n}{0}{1}, \label{eqn:Aut_Cn} \end{equation}
and elementary arithmetic considerations, which we omit here for brevity, then show that we have a bijection
\begin{equation} \begin{array}{ccc} \{ (c,d) \in \Z/N\Z \times (\Z/(c,N)\Z)^\times \} / \pm 1 & \longleftrightarrow  & \text{Cusps}(X_1(N)) \\ (c,d) & \longmapsto & \big( C_{N/(c,N)}, (\zeta_N^d, c/(c,N) ) \big), \end{array} \label{eqn:cusps_X1} \end{equation}
where for brevity we have written~$(c,N)$ for~$\gcd(c,N)$, and where~$\zeta_N$ is a fixed primitive~$N$-th root of unity.
We may thus represent the cusps of~$X_1(N)$ by such pairs~$(c,d)$ up to negation.

The advantage of this representation is that the Galois action on the cusps is then transparent. Indeed,~\eqref{eqn:cusps_X1} confirms that the cusps are all defined over~$\Q(\mu_N)$, and shows that for each~$x \in (\Z/N\Z)^\times$, we have
\begin{equation} \sigma_x \cdot (c,d) = (c,xd) \label{eqn:cusps_Gal} \end{equation}
where~$\sigma_x \in \Gal(\Q(\mu_N)/\Q)$ is as in~\eqref{eqn:Gal_cyclo}. In particular, two cusps~$(c,d)$,~$(c',d')$ of~$X_1(N)$ are in the same Galois orbit iff.~$c = \pm c' \bmod N$.

One easily verifies that the correspondence between this representation of the cusps by pairs~$(c,d)$ and the more traditional one by classes of elements of~$\Q \cup \{ \infty \}$ is as follows: given an element~$a/c \in \Q \cup \{ \infty \}$ in lowest terms, find~$b,d \in \Z$ such that~$\smatabcd \in \SL_2(\Z)$; then the cusp represented by~$a/c$ in the traditional representation is represented by~$(c,d)$ in our representation, and vice-versa.

\begin{ex}\label{ex:cusps_X1}
The cusp~$\infty = 1/0$ is represented by~$(c=0,d=1)$. In particular,~\eqref{eqn:cusps_Gal} shows that this cusp is fixed by~$\sigma_x$ only if~$x=\pm1$, which means that its field of definition is~$\Q(\zeta_N+\zeta_N^{-1})$. This can be visualised by noticing that it corresponds under~\eqref{eqn:cusps_X1} to the pair~$\big(C_1, (\zeta_N,0)\big)$, and by applying~\eqref{eqn:Aut_Cn} to~$n=N/(c,N)=1$, which shows that this pair is isomorphic (by negation) to~$\big(C_1, (\zeta_N^{-1},0)\big)$ but not to any other of its Galois conjugates.

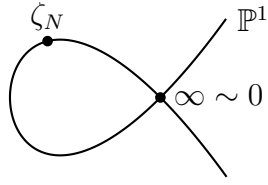
\begin{figure}[H]
\begin{center}
\begin{tikzpicture}[scale=2]
\draw[thick,variable=\t,domain=-1.2:1.2,samples=500]
  plot ({\t*\t - 1 )},{\t*\t*\t - \t})
  node[right] {$\PP^1$};
\fill[black] (0,0) circle (1pt)
  node [right, xshift=1pt, yshift=1pt] {$\infty \sim 0$};
\fill[black] (-0.75,0.375) circle (1pt)
  node [above] {$\zeta_N$};
\end{tikzpicture}
\end{center}
\caption*{Moduli interpretation of the cusp~$\infty$}
\end{figure}
We may interpret this by noticing that for each~$\tau$ in the upper half-plane, we have a point on~$X(N)(\C)$ corresponding to the pair~$(\C/\langle \tau,1 \rangle, \beta)$, where~$\beta(1,0) = \tau/N$ and~$\beta(0,1)=1/N$ so that~$e_N(\beta(1,0),\beta(0,1))=e^{2\pi i /N}$ by our normalisation~\eqref{eqn:normalisation_Weil} of the Weil pairing. This pair projects by~\eqref{eqn:proj_XN_XH} to the point of~$X_1(N)$ represented by~$(\C/\langle \tau,1 \rangle, 1/N) \simeq (\C^\times/q^\Z,e^{2\pi i /N})$ where~$q=e^{2\pi i \tau}$, so when~$\tau \rightarrow \infty$, it becomes~$(\G_m,e^{2\pi i /N})$, which is not defined over~$\Q$.

\bigskip

On the contrary, the cusp~$0=0/1$ is represented by~$(c=1,d=0)$, and is thus defined over~$\Q$; indeed, it corresponds to the pair~$\big(C_N, (1,1)\big)$, which is clearly defined over~$\Q$.

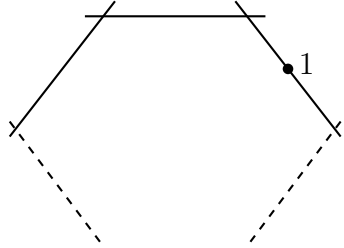
\begin{figure}[H]
\begin{center}
\begin{tikzpicture}[domain=-2:2]
\draw[thick] (-1.2,1) -- (1.2,1);
\draw[thick] (-0.8,1.2) -- (-2.2,-0.6);
\draw[thick,dashed] (-2.2,-0.4) -- (-1,-2);
\draw[thick] (0.8,1.2) -- (2.2,-0.6);
\fill[black] (1.5,0.3) circle (2pt)
  node [right,yshift=2pt] {$1$};
\draw[thick,dashed] (2.2,-0.4) -- (1,-2);
\end{tikzpicture}
\end{center}
\caption*{Moduli interpretation of the cusp~$0$}
\end{figure}

\end{ex}

As explained in the previous subsection, in this article, we will actually not work with~$X_1(N)$, but rather with~$X_H(N)$ where~$H$ is a subgroup of~$(\Z/N\Z)^\times$ containing~$-1$. Fortunately, translating the above results to the case of the cusps of~$X_H(N)$ presents no difficulty. Indeed, going down from~$X_1(N)$ to~$X_H(N)$ amounts to identifying
\begin{equation} \gamma \cdot s \sim s \label{eqn:ident_cusps_X1_XH} \end{equation}
for each cusp~$s$ of~$X_1(N)$ and each~$\gamma \in \SL_2(\Z)$ congruent mod~$N$ to~$\smat{h^{-1}}{*}{0}{h}$ for some~$h \in H$. Given such a cusp~$s$ represented by~$(c,d)$, let~$M_s =\smat{*}{*}{c}{d} \in \SLZ$ be such that~$M_s \cdot \infty = s$; then, given such a~$\gamma$, we find
\begin{equation} \gamma \mat{*}{*}{c}{d} \equiv \mat{*}{*}{hc}{hd} \bmod N. \label{eqn:Gamma0_on_cusps} \end{equation}
This means that under our representation of cusps by pairs~$(c,d)$,~\eqref{eqn:ident_cusps_X1_XH} becomes
\[ (hc,hd) \sim (c,d) \ \forall h \in H. \]
Therefore,~\eqref{eqn:cusps_X1} simply becomes
\begin{equation} \{ (c,d) \in \Z/N\Z \times (\Z/(c,N)\Z)^\times \} / H \longleftrightarrow \text{Cusps}(X_H(N)); \label{eqn:cusps_XH} \end{equation}
in others words, we now consider the pairs~$(c,d)$ up to multiplication by~$H$ instead of up to negation. The same computation also shows that for all $y \in (\Z/N\Z)^\times$, the diamond operator~$\langle y \rangle$ takes the cusp represented by~$(c,d)$ to that represented by~$(yc,yd)$.

\begin{ex}\label{ex:cusps_XH}
Continuing example~\ref{ex:cusps_X1}, we see that the field of definition of the cusp~$\infty$ of~$X_H(N)$ is~$\Q(\mu_N)^H$, where we view~$H$ as a subgroup of~$\Gal\big(\Q(\mu_N)/\Q\big)$ thanks to~\eqref{eqn:Gal_cyclo}.
\end{ex}

\begin{rk}\label{rk:cusps_XH_rat}
The description~\eqref{eqn:cusps_XH} of the cusps and~\eqref{eqn:cusps_Gal} of the Galois action on them shows that the modular curves~$X_H(N)$ tend to have a large supply of rational cusps.
\end{rk}

\subsubsection*{Widths}

In order to be able to consider~$q$-expansions at various cusps, we must determine the width of these cusps.

Let~$s$ be a cusp of~$X_H(N)$, and let again~$M_s \in \SLZ$ be such that~$M_s \cdot \infty = s$. The width of~$s$ is then by definition the smallest positive integer~$w$ such that
\[ M_s \smat{1}{w}{0}{1} M_s^{-1} \in \Gamma_H(N). \]
Writing~$M_s = \smatabcd$, we compute
\[ M_s \mat{1}{w}{0}{1} M_s^{-1} = \mat{1-acw}{a^2w}{-c^2w}{1+acw}, \]
so we want~$c^2 w \equiv 0 \bmod N$ and~$1 \pm acw \in H$ (note that~$(1+acw)(1-acw)=1-a^2c^2w^2 \equiv 1$ if~$c^2w \equiv 0$, so the~$\pm$ sign is irrelevant, i.e. this identity with either sign implies the one with the other sign). Let~$g_2 = \gcd(N,c^2)$ and~$N_2 = N/g_2$; then~$c^2 w \equiv 0$ iff.~$N_2 \mid w$, whence finally
\begin{equation} w = N_2 \min \{ t \in \N \ \vert \ 1+acN_2 t \in H \}. \label{eqn:width} \end{equation}

This formula allows us to determine the width of a cusp represented by a pair~$(c,d)$. 
Unfortunately, it is a bit tedious to apply for general~$H$, but it simplifies considerably if work with~$X_0(N)$ (which amounts to~$H=(\Z/N\Z)^\times$) or with~$X_1(N)$ (which amounts to~$H = \{ \pm 1 \}$). For future reference, we note the following result, which is valid for general~$H$:

\begin{pro}\label{pro:w=N}
The cusp represented by~$(c,d)$ has width~$w=N$ iff.~$\gcd(c,N)=1$.
\end{pro}

\begin{proof} Let~$w$ be the with of the cusp represented by~$(c,d)$, and let~$g=\gcd(c,N)$, so that~$ g \mid g_2 \mid g^2$.

Suppose first that~$g>1$; then~$g_2>1$ so~$N_2<N$. We distinguish two cases: if~$g < g_2$, then taking~$t=g$ yields
\[ 1+acN_2t = 1+a \frac{c}g N \frac{g^2}{g_2} = 1, \]
so the smallest possible~$t$ is at most~$g$ whence~$w \leqslant N_2 g < N_2 g_2 = N$; and if~$g=g_2$, then taking~$t=1$ we get
\[ 1+acN_2t = 1+a \frac{c}g N = 1 \in H, \]
so the smallest possible~$t$ is~$t=1$ whence~$w=N_2 < N$.

Conversely, if~$g=1$, then~$g_2=1$, so~$N_2=N$ whence~$w=N$.
\end{proof}

\subsection{Rationality of~$q$-expansions}\label{sect:Qqexp}

In order to compute modular Galois representations, we will need to construct rational maps~$J_H(N) \dashrightarrow \A^1$. As we will see in section~\ref{sect:Eval}, one way to do so involves looking at the~$q$-expansion coefficients~$a_n(f)$ of some forms~$f$ at some cusp; but it is of course fundamental for our purpose that the dependency of the~$a_n(f)$ on~$f$ be Galois-equivariant. This is unfortunately not the case at every cusp; for instance, it is not the case of the cusp~$\infty$ of~$X_1(N)$ for~$N > 6$ since this cusp is not even defined over~$\Q$.

In fact, thinking about the rationality of the~$q$-expansion in terms of the cusp alone is wrong. Indeed, given a function~$f \in \overline \Q\big(X_H(N)\big)$ (or more generally, a modular form) and a cusp~$s$ of~$X_H(N)$ of width~$w$, it is tempting to define ``the''~$q$-expansion of~$f$ at~$s$ as the expansion of~$f \mact M_s$ at~$\infty$ in terms of~$q_w = e^{2 \pi i \tau/w}$, where~$M_s \in \SLZ$ satisfies~$M_s \cdot \infty = s$. However, this definition does not make actual sense if~$w>1$. Indeed, the matrices~$M \in \SLZ$ satisfying~$M \cdot \infty = s$ are precisely those of the form~$M = \pm M_s \smat{1}{x}{0}{1}$ where~$x \in \Z$, and while the~$\pm$ sign does not matter, different values of~$x$ yield different~$q_w$-expansions of~$f \mact M$ at~$\infty$; more precisely, we have
\begin{equation} a_n\left(f \mact M_s \smat{1}{x}{0}{1}\right) = e^{2 \pi i n /w} a_n\left(f \mact M_s\right) \label{eqn:loc_param} \end{equation}
for all~$n \in \Z$. However, this still shows that the coefficient~$a_0(f \mact M_s)$ and the order of vanishing of~$f$ do not depend on~$M_s$, but only on~$s$.

We are thus led to the following definition:

\begin{de}\label{de:Qqexp}
We say that a matrix~$M \in \SLZN$  \emph{yields rational~$q$-expansions} if the map
\[ f \longmapsto q\text{-expansion of } (f\mact M) \text{ at } \infty \]
is Galois-equivariant.
\end{de}

This condition is equivalent to the requirement that the~$q$-expansion of~$f \mact M$ have rational coefficients whenever~$f~$ is defined over~$\Q$.

\begin{rk}
Technically, these expansions are~$q_w$-expansions, where~$w \in N$ is the width of the cusp~$M \cdot \infty$. One way to circumvent this technicality would be to talk about~$q_N$-expansions, since~$w \mid N$ always. However, this does not impact our discussion about the Galois-equivariance of the coefficients, so for convenience, we will persist in this abuse of language in the rest of this section.
\end{rk}

We must therefore determine which~$M \in \SLZN$ yield rational~$q$-expansions. Recall from~\cite[6.2]{Shimura} that every element~$f \in F_N$ has a (possibly Laurent)~$q_N$-expansion
\[ f = \sum_{n \gg -\infty} a_n(f) q_N^n \]
with coefficients~$a_n(f) \in \Q(\mu_N)$, and that we have the relation
\begin{equation} a_n(f)^{\sigma_x} = a_n(f \mact \smat{1}{0}{0}{x}) \label{eqn:Gal_cyclo_qexp} \end{equation}
for all~$x \in (\Z/N\Z)^\times$ and~$n \in \Z$, where~$\sigma_x \in \Gal(\Q(\mu_N)/\Q)$ is as in~\eqref{eqn:Gal_cyclo}.

From~\eqref{eqn:fnfield_XH}, we deduce that for all~$M \in \SLZN$,
\begin{align*}  & M \text{ yields rational~$q$-expansions} \\
\Longleftrightarrow \ & \forall f \in \Q\big( X_H(N)\big) , \ \forall n \in \Z, \ a_n(f \mact M) \in \Q \\
\Longleftrightarrow \ & \forall f \in \Q\big( X_H(N)\big), \ \forall x \in (\Z/N\Z)^\times, \ f \mact M = f \mact M \smat{1}{0}{0}{x} \\
\Longleftrightarrow \ & \forall x \in (\Z/N\Z)^\times, M \smat{1}{0}{0}{x} M^{-1} \in U_H. \\
\end{align*}
Writing~$M=\matabcd$, this translates explicitly into
\begin{equation} \forall x \in (\Z/N\Z)^\times, \ cd(x-1) = 0 \bmod N \text{ and } ad(x-1)+1 \in H. \label{eqn:crit_M_rat_qexp} \end{equation}
This criterion allows us to determine explicitly for which cusps~$s$ of~$X_H(N)$ there exists~$M \in \SLZN$ such that~$M \cdot \infty = s$ and that~$M$ yields rational~$q$-expansions, and to find such an~$M$ if one exists. 

\begin{rk}\label{rk:always_1_Qcusp}
There is always at least one such cusp; namely, we can take~$s=0$ and~$M = \smat{0}{-1}{1}{0}$, since~\eqref{eqn:crit_M_rat_qexp} is obviously satisfied if~$d=0$.
\end{rk}

\bigskip

Similarly to formula~\eqref{eqn:width}, criterion~\eqref{eqn:crit_M_rat_qexp} is a bit tedious to use in practice for general~$H$, but simplifies considerably if we work with~$X_0(N)$ or~$X_1(N)$. For instance, we have the following results:

\begin{pro}\label{pro:ratq_1}
Suppose~$H=\{\pm1\}$ and~$\phi(N) \geqslant 3$. Then~$M$ yields rational~$q$-expansions iff.~$2d = 0 \bmod N$.
\end{pro}

\begin{proof}
In this case, we have~$X_H(N) = X_1(N)$, and~$\Q\big(X_1(N)\big) = \Q(j,f^{(0,1)}_0)$ according to~\cite[7.7]{DS}, where the~$f_0^v$ were defined at the beginning of this section; therefore~$M = \smatabcd \in \SLZN$ yields rational~$q$-expansions iff.~$f^{(0,1)}_0 \mact M \smat{1}{0}{0}{x} = f^{(0,1)}_0\mact M$ for all~$x \in (\Z/N\Z)^\times$, which translates into 
\begin{equation} (c,dx) \equiv \pm (c,d) \text{ for all } x \in (\Z/N\Z)^\times. \label{eqn:crit_M_rat_qexp_1} \end{equation}
We now distinguish three cases.

If~$c \not \equiv -c$, then for all~$x$,~$dx \equiv d$, i.e.~$N \mid d(x-1)$. In particular, taking~$x=-1$ shows that~$N \mid 2d$. Conversely, assume that~$N \mid 2d$. If~$N$ is odd, then~$N \mid d$, so~\eqref{eqn:crit_M_rat_qexp_1} is clearly satisfied; and if~$N$ is even, then~$N/2 \mid d$, so~$N \mid d(x-1)$ since~$(x-1)$ is even for all~$x \in (\Z/N\Z)^\times$.

If~$c \equiv -c$ but~$c \not \equiv 0$, then~$N$ is even and~$c \equiv N/2$, so that~$1 = \gcd(c,d,N) = \gcd(d,N/2)$. Let~$x \in (\Z/N\Z)^\times$; then~$x+1$ and~$x-1$ are even, so~\eqref{eqn:crit_M_rat_qexp_1} implies implies~$N \mid d(x \pm 1)$ whence~$\frac{N}2 \mid d \frac{x \pm 1}2$ so~$\frac{N}2 \mid \frac{x \pm 1}2$ so~$N \mid (x \pm 1)$ so~$x \equiv \pm 1$, which contradicts our assumption that~$\phi(N) \geqslant 3$. So this case cannot happen.

Finally, if~$c \equiv 0$, then~$1=\gcd(d,N)$, and~\eqref{eqn:crit_M_rat_qexp_1} implies that for all~$x \in (\Z/N\Z)^\times$,~$d(x+1)$ or~$d(x-1)$ vanishes mod~$N$. Since~$d$ is invertible mod~$N$, this implies~$x \equiv \pm 1$, so again this case cannot occur.
\end{proof}

\begin{cor}\label{cor:rat_w=N}
Let~$s$ be a cusp of~$X_1(N)$. Suppose that~$\phi(N) \geqslant 3$, and that~$N$ is either odd or a multiple of~$4$. Then there exists~$M \in \SLZN$ such that~$M \cdot \infty = s$ and that~$M$ yields rational~$q$-expansions iff.~$s$ has width~$N$.
\end{cor}

\begin{proof}
Let~$(c,d)$ represent the cusp~$s$. Since~$d$ lives in~$(\Z/(c,N)\Z)^\times$, we may replace it with~$d_t=d+tc$ for any~$t \in \Z$.

Suppose~$s$ has width~$N$. Then~$\gcd(c,N)=1$ by proposition~\ref{pro:w=N}. Let~$(c,d)$ represent the cusp~$s$. Since~$d$ lives in~$(\Z/(c,N)\Z)^\times$, we may replace it with~$d_t=d+tc$ for any~$t \in \Z$. As~$\gcd(c,N)=1$, we can choose~$t$ so that~$d_t \equiv 0 \bmod N$, and then~\eqref{eqn:crit_M_rat_qexp} obviously holds.

Conversely, suppose there is such an~$M = \smatabcd \in \SLZN$, so that~$s$ is represented by~$(c,d)$ and that~$N \mid 2 d$. We are going to deduce that~$\gcd(c,N)=1$, which will conclude by proposition~\ref{pro:w=N}. Observe that~$1 = \gcd(c,d,N)$.  We distinguish two cases: if~$N$ is odd, then~$N \mid d$, so~$\gcd(c,N)=1$; and if~$4 \mid N$, then~$\frac{N}2 \mid d$, so that~$\gcd(c,N/2)$ divides~$\gcd(c,d,N ) = 1$ and is therefore~$1$, and since~$N/2$ is even, this implies~$\gcd(c,N)=1$.
\end{proof}

\begin{rk}
Corollary~\ref{cor:rat_w=N} may fail if~$2 \mid N$ but~$4 \nmid N$, as illustrated by the example~$N=10$,~$s$ represented by~$(c=2,d=5)$. Indeed, this cusp has width~$w=5$ by~\eqref{eqn:width}, and yet the matrix~$M = \smat{1}{2}{2}{5} \in \SLZN$ yields rational~$q$-expansions by~\eqref{eqn:crit_M_rat_qexp} and takes~$\infty$ to~$s$. 
\end{rk}

\section{Makdisi's moduli-friendly Eisenstein series}\label{sect:MakEis}

\subsection{Makdisi's construction}

In order to construct~$p$-adic Makdisi models of modular curves without resorting to explicit plane models, we will rely on ``moduli-friendly'' modular forms in the sense of~\cite{MakEis}, meaning that their ``value'' at a point of the modular curve can be easily read off the representation of this point as an elliptic curve equipped with some appropriate level structure. We thus think of our modular forms ``\`a la Katz''; in other words, we view a modular form~$f$ of weight~$k$ and level~$\Gamma(N)$ over a ring~$R$ in which~$N$ is invertible as a function on the set of isomorphism classes of triples
\[ (E,\omega,\beta) \]
where~$E$ is an elliptic curve over~$R$,~$\omega$ is a generator of the sheaf of regular relative differentials on~$E/R$, and~$\beta$ is what~\cite[2.0.3]{Katz} calls a \emph{na\"ive level~$N$ structure} on~$E$, that is to say an isomorphism
\[ \beta : (\Z/N\Z)^2 \simeq E[N] \]
of group schemes over~$R$, and satisfying the homogeneity condition
\[ f(E,\lambda \omega,\beta) = \lambda^{-k} f(E,\omega,\beta) \]
for all~$\lambda \in R^\times$ as well as some extra compatibility conditions (namely commutation with base change, cf.~\cite[2.1]{Katz} for details).

We will in fact restrict ourselves to elliptic curves defined by short Weierstrass equations
\begin{equation} (\W) : y^2 = x^3+Ax+B. \label{eqn:WeiW} \end{equation}
By assigning to such an equation the differential~$\omega_{\W} = dx/2y$, it then makes sense to talk about the ``value''~$f(\W,\beta)$ of~$f$ at the pair~$(\W,\beta)$; in other words, choosing a short Weierstrass model for~$E$ yields a local trivialisation of the sheaf of modular forms.

\bigskip

Fix a level~$N \in \N$, and let~$R$ be a ring in which~$6N$ is invertible. In~\cite{MakEis}, Makdisi constructs Eisenstein forms~$f_1^v$ of weight~$k=1$ and level~$\Gamma(N)$ over~$R$ indexed by non-zero vectors~$v \in (\Z/N\Z)^2$ and which enjoy the following particularly nice properties:

\begin{thm}\label{thm:MakEis}\mbox{}
\begin{itemize}
\item (Moduli-friendliness) Let~$v_1, v_2 \in (\Z/N\Z)^2$ be such that neither~$v_1$, nor~$v_2$, nor~$v_3=-v_1-v_2$ are zero. Given a pair~$(\W,\beta)$ where~$\W$ is a short Weierstrass equation defining an elliptic curve~$E$ over~$R$ and~$\beta$ is a na\"ive level~$N$ structure on~$E$, the ``value''
\[ f_1^{v_1}(\W,\beta) + f_1^{v_2}(\W,\beta) + f_1^{v_3}(\W,\beta) \]
agrees with the slope of the line joining the aligned points~$\beta(v_1)$,~$\beta(v_2)$, and~$\beta(v_3)$ on the model of~$E$ defined by~$\W$ (to be interpreted as the slope of the flex tangent in the case where~$v_1=v_2=-v_1-v_2$).

\item (Generation) The subalgebra of the~$R$-algebra
\[ \bigoplus_{k \geqslant 0} \M_k\big(\Gamma(N);R\big) \]
generated by the~$f_1^v$ is
\[ R \oplus \Ei_1\big(\Gamma(N),R\big) \oplus  \bigoplus_{k \geqslant 2} \M_k\big(\Gamma(N);R\big); \]
in other words, as~\cite{MakEis} puts it, it ``misses'' precisely the cusp forms of weight 1.
\end{itemize}
\end{thm}

Whenever~$P,Q \in E$ are points such that
\begin{equation} \text{neither~$P$, nor~$Q$, nor~$P+Q$ are at infinity,} \label{cond:slope} \end{equation}
denote by~$\lambda_{P,Q}$ the slope of the line joining~$P$ and~$Q$ on the model~$\W$ of~$E$ if~$P\neq Q$, and the slope of the tangent line of~$E$ at~$P$ if~$P=Q$; observe that~\eqref{cond:slope} ensures that this line is not vertical, so that this slope is well defined. The first property can be summarised by
\begin{equation} f_1^{v_1}(\W,\beta) + f_1^{v_2}(\W,\beta) + f_1^{-v_1-v_2}(\W,\beta) = \lambda_{\beta(v_1),\beta(v_2)} \label{eqn:l1_l2} \end{equation}
even in the case where~$v_1$,~$v_2$ and~$-v_1-v_2$ are not distinct. Makdisi shows that this relation can be inverted so as to read the ``value'' of the~$f_1^v$ off the slopes of the lines joining the~$N$-torsion points of~$E$, for instance as
\begin{equation} f_1^v(\W,\beta) = \frac1N \sum_{x \bmod N} \lambda_{\beta(v),\beta(v+xw)} \label{eqn:f1_inversion} \end{equation}
where~$w$ is any vector of~$(\Z/N\Z)^2$ whose span intersects trivially that of~$v$, so that the slope~$\lambda_{\beta(v),\beta(v+xw)}$ is well-defined for all~$x$ (however, see below for a more efficient method). Thus the~$f_1^v$ are truly moduli-friendly modular forms.

By combining this observation with the second property, we thus get the neat statement that apart from cusp forms of weight 1, any modular form can (in principle) be expressed as a polynomial in the moduli-friendly forms~$f_1^v$, and therefore evaluated at a pair~$(\W,\beta)$.

\subsection{Efficient evaluation of $f_1^v$}

In order to construct~$p$-adic Makdisi models of modular Jacobians, we will use these modular forms over rings of the form~$R=\Z_q/p^e$ where~$p \nmid 6N$, in which the most computationally expensive of the four basic operations is by far division. Using~\eqref{eqn:f1_inversion} to evaluate such a form~$f_1^v$ at a pair~$(\W,\beta)$ requires us to determine~$N$ slopes of lines joining~$N$-torsion points on~$E$, or of tangent lines at such points. Evaluating such a slope requires one division in~$R$, even in the case of a tangent line, since differentiating~\eqref{eqn:WeiW} shows that the slope of the tangent to~$E$ at the point~$P=(x,y)$ is~$\lambda_{P,P} = \frac{3x^2+A}{2y}$. In total, each evaluation of~$f_1^v$ by~\eqref{eqn:f1_inversion} therefore requires~$N$ divisions in~$R$, which may be prohibitively costly if~$e$ and~$N$ are large. We will therefore use a better approach, which was hinted at by Makdisi in \cite[3.14]{MakEis} and allows one to evaluate~$f_1^v$ by evaluating only~$O(\log N)$ slopes. We now describe this approach in detail, basing ourselves on explanations provided by Makdiski to the author.

\begin{pro}\label{pro:Kamal_log}
Let~$(\W,\beta)$ be a pair as above, let~$v \in (\Z/N\Z)^2$ be a non-zero vector, and let~$n | N$ be the exact order of~$v$, so that~$P = \beta(v)$ is a point of~$E[N]$ of exact order~$n$. Let~$c_1,c_2,\cdots,c_{n-1}$ be the finite~$R$-valued sequence defined by
\begin{equation} c_m = \left\{ \begin{array}{ll}
0 & \text{ if } m=1,\\
2 c_{m/2} + \lambda_{P_{m/2},P_{m/2}} & \text{ if } 1<m<n \text{ is even},\\
c_{m-1} + \lambda_{P,P_{m-1}} & \text{ if } 1<m<n \text{ is odd},\\
\end{array} \right. \label{eqn:Kamal_log} \end{equation}
where for brevity we have written~$P_m$ for~$[m]P = \beta(mv) \in E[N]$. Then
\begin{itemize}
\item[(i)] Condition~\eqref{cond:slope} is satisfied in every case of~\eqref{eqn:Kamal_log}, so that this sequence is well-defined,
\item[(ii)] \label{pro:Kamal_log:ind} For all~$1 \leqslant m < n$, we have~$f_1^{mv}(\W,\beta) = m f_1^{v}(\W,\beta) - c_m$,
\item[(iii)]~$\displaystyle f_1^{v}(\W,\beta) = \frac1n c_{n-1}$.
\end{itemize}
\end{pro}

\begin{proof}
\begin{itemize}
\item[(i)] In the case where~$1 < m < n$ is even, neither~$P_{m/2}$ nor~$P_{m/2}+P_{m/2}$ are at infinity since~$m/2<m<n$. In the case where~$1 < m < n$ is odd, neither~$P$ nor~$P_{m-1}$ nor~$P+P_{m-1} = P_m$ are at infinity since~$1 \leqslant m-1 < m < n$.
\item[(ii)] First of all, observe that~$f_1^v$ is an odd function of~$v$, meaning that~$f_1^{-v} = -f_1^v$ for all~$v$; this is apparent on~\eqref{eqn:f1_inversion}, and actually follows directly from Makdisi's construction. The formula~$f_1^{mv}(\W,\beta) = m f_1^{v}(\W,\beta) - c_m$ then follows by induction on~$m<n$. Indeed, it obviously holds for~$m=1$. Suppose now that it holds for all~$m'<m$. In the case where~$m$ is even,~\eqref{eqn:l1_l2} applied with~$v_1 = v_2 = \frac{m}2 v$ shows that 
\[ 2 f_1^{\frac{m}2 v}(\W,\beta) - f_1^{mv}(\W,\beta) = \lambda_{P_{m/2},P_{m/2}}, \]
whence~$f_1^{mv}(\W,\beta) = 2 \frac{m}2 (f_1^v(\W,\beta)-c_{m/2}) - \lambda_{P_{m/2},P_{m/2}}$. Similarly, in the case when~$m$ is odd,~\eqref{eqn:l1_l2} applied with~$v_1=v$,~$v_2=(m-1)v$ yields
\[ f_1^v(\W,\beta) + f_1^{(m-1)v}(\W,\beta)-f_1^{mv}(\W,\beta) = \lambda_{P,P_{m-1}}, \]
whence~$f_1^{mv}(\W,\beta) = f_1^v(\W,\beta) + \big((m-1) f_1^v(\W,\beta) - c_{m-1}) - \lambda_{P,P_{m-1}}$.
\item[(iii)] Taking~$m=n-1$ in (ii) and using again the fact that~$f_1^v$ is an odd function of~$v$, we obtain
\[ -f_1^v(\W,\beta) = f_1^{(n-1)v} (\W,\beta) = (n-1) f^1_v(\W,\beta) - c_{n-1}. \qedhere \]
\end{itemize}
\end{proof}

\begin{rk}\label{rk:l1_l2_1}
The formula~\eqref{eqn:f1_inversion} and the algorithm outlined in proposition~\ref{pro:Kamal_log} both demonstrate that the forms~$l_1^{v,w} : (\W,\beta) \mapsto \lambda_{\beta(v),\beta(w)}$ span the same~$R$-algebra of modular forms as the~$f_1^v$. However, although the generators~$l_1^{v,w}$ may seem more appealing since they are easier to evaluate than the~$f_1^v$, we shall demonstrate in remark~\eqref{rk:l1_l2_2} below that using the~$f_1^v$ results in a better complexity in the construction of ~$p$-adic Makdisi models of modular curves, whence our focus on the~$f_1^v$ in this section.
\end{rk}

\section{Makdisi models of modular Jacobians}\label{sect:MakMod}

\subsection{Strategy}\label{sect:Mak_model_strategy}

We now have all the ingredients required to construct~$p$-adic Makdisi models of modular curves. Since we want to compute modular Galois representations, as explained in section~\ref{sect:rho_in_XH} we focus on the case of the curves~$X_H(N)$, where~$N \in \N$ and~$H \ni -1$ is a subgroup of~$(\Z/N\Z)^\times$.

For simplicity, we make the following assumptions:
\begin{itemize}
\item~$X_H(N)$ has at least 3 cusps,
\item~$p$ does not divide~$6$, nor~$\ell$, nor~$N$, nor the order of the subgroup~$H \leqslant (\Z/N\Z)^\times$.
\end{itemize}

We will explain the reason for these assumptions below; for now, we just note that the assumption on the number on cusps is not an essential one (cf. remark~\ref{rk:3cusps} below), and that we require~$p \nmid N$ since~$p$-adic Makdisi models require~$p$ to be a prime of good reduction,~$p \neq \ell$ since our method relies on the~$\ell$-torsion being \'etale at~$p$, and that~$p \nmid 6$ so as to ensure the validity of Makdisi's construction of the moduli-friendly Eisenstein series~$f_1^v$. We will explain how the prime $p$ is chosen in section~\ref{sect:Choice_p} below.

\bigskip

We can determine the local L factor of~$X_H(N)$ at~$p$ by~\eqref{eqn:Lp_mod_res}, from which we recover in particular the the genus~$g \in \N$ of~$X_H(N)$.

\bigskip

In order to construct a Makdisi~$p$-adic model, we then need to pick a line bundle~$\L$ on~$X_H(N)$. It is natural to choose a line bundle whose sections are modular forms of level~$\Gamma_H(N)$; we choose~$\L$ so that its sections are the modular forms (not just cusp forms) of weight 2. This means that the degree of~$\L$ is~$2g-2+\nu_\infty$, where~$\nu_\infty$ is the number of cusps of~$X_H(N)$; the assumption~$\nu_\infty \geqslant 3$ that we have made above thus ensures that the requirement~\eqref{eqn:d0bound} is met.

\begin{rk}\label{rk:3cusps}
If our modular curve happens to have fewer than~3 cusps, we can still apply the same construction, by choosing~$\L$ so that its sections are modular forms of some higher weight, thus ensuring that~$\deg \L$ is large enough; besides, the optimisation process presented in subection~\ref{sect:pruning} will still apply up to straightforward modifications since there is always at least one cusp. However, in almost all the cases that will be relevant to us in this article,~$\nu_\infty$ is much larger than~$3$, so we make the assumption that~$\nu_\infty \geqslant 3$ for the simplicity of the exposition.
\end{rk}

As in the previous section, we view the non-cuspidal points of~$X_H(N)$ as pairs~$(\W,H \cdot P)$, where~$\W$ is a Weierstrass equation defining an elliptic curve~$E$, and~$P \in E[N]$. In order to construct our~$p$-adic Makdisi model, we need to fix sufficiently many such points at which to evaluate (under some local trivialisation of~$\L$) a basis of global sections of~$\L$. Choosing different elliptic curves~$E$ would require keeping track of the~$N$-torsion point~$P$ as one elliptic curve deforms into another, which seems complicated. Instead, Makdisi brilliantly suggests to \emph{fix} the curve~$E$, and to consider the points on the modular curve corresponding to the various possible~$N$-torsion points~$P$ on that~$E$; in other words, to work in a fibre of the projection map~$\pi : X_H(N) \longrightarrow X(1)$.

For simplicity, we choose~$\W$ so that it (not just~$E$) has good reduction at~$p$. Since~$p \nmid N$ by assumption, N\'eron-Ogg-Shafarevic ensures that the coordinates of the~$N$-torsion points of~$E$ generate an unramified extension~$\Q_q$ of~$\Q_p$. 
Furthermore, our assumption that~$\W$ has good reduction at~$p$ ensures that the coordinates of the nonzero~$N$-torsion points of~$E$ actually lie in the ring of integers~$\Z_q$ of~$\Q_q$.

\bigskip

We thus fix~$A, B \in \Z$ defining an elliptic curve
\[ (\W) : y^2 = x^3+Ax+B \]
over~$\Q$ having good reduction at~$p$ and whose~$j$-invariant is neither~$0$ nor~$1728$ mod~$p$, so as to avoid the ramification locus of~$\pi$. The number of points in the fibre of~$\pi$ above~$E$ is then equal to the degree~$d$ of~$\pi$. By~\cite[3.1.1]{DS}, the genus of~$X_H(N)$ is
\[ g = 1 + \frac{1}{12} d - \frac14{\nu_2} - \frac1{\nu_3} - \frac12{\nu_\infty} \]
where~$\nu_2$ (resp.~$\nu_3$) denotes the number of elliptic points of~$X_H(N)$ of order~$2$ (resp.~$3$). It follows that
\[ d_0 = 2g-2+\nu_\infty = \frac16 d - \frac12 \nu_2 - \frac23 \nu_3, \]
whence
\[ d-5d_0 = \frac16 d + \frac52 \nu_2 + \frac{10}3 \nu_3 > 0, \]
which shows that the lower bound~\eqref{eqn:nZbound} on the number of points of~$X_H(N)$ at which we evaluate the sections of~$\L$ is satisfied; we are thus in good shape to construct a valid Makdisi~$p$-adic model of~$X_H(N)$.

\bigskip

We then determine the coordinates in~$\F_q$ of the~$N$-torsion points of~$E$ in the model~$\W$, and, having set a desired accuracy~$O(p^e)$ for our~$p$-adic Makdisi model, we Hensel-lift these coordinates to~$\Z_q/p^e$. Besides, we arbitrarily fix a level structure~$\beta : (\Z/N\Z)^2 \simeq E[N]$. By~\eqref{eqn:fibre_XH}, the points at which we evaluate our forms, that is to say the points on the fibre of the projection~$X_H(N) \rightarrow X(1)$ at~$E$, may then be identified with the primitive vectors of~$(\Z/N\Z)^2$ up to scaling by~$H$. In particular, if~$\Phi \in \GLZN$ is the matrix describing the action of the Frobenius~$\Frob_p$ on~$E[N]$ with respect to~$\beta$, then the image by~$\Frob_p$ of the point of the fibre corresponding to the primitive vector~$v \in (\Z/N\Z)^2$ is the point of the fibre corresponding to the vector~$v \cdot (^t \Phi)$. We can therefore determine the permutation induced by~$\Frob_p$ on the fibre, provided that we have computed the matrix~$\Phi$. We explain in detail how all this is done in subsection~\ref{sect:LiftEN} below.

\bigskip

By the second part of theorem~\ref{thm:MakEis}, the space of modular forms of weight 2 and level~$\Gamma(N)$ is spanned by the products~$f_1^v f_1^w$, where~$v$ and~$w$ range over the set of nonzero vectors of~$(\Z/N\Z)^2$. We use these products~$f_1^v f_1^w$ to construct weight-2 forms of level~$\Gamma_H(N)$ by taking traces; namely, we set
\begin{equation} f_{2,H}^{v,w} = \Tr^{\Gamma(N)}_{\Gamma_H(N)} f_1^v f_1^w \overset{\text{def}}{=} \sum_{\gamma \in \overline \Gamma_H(N)} (f_1^v f_1^w) \mact \gamma =  \sum_{\gamma \in \overline \Gamma_H(N)} f_1^{v \gamma} f_1^{w \gamma}, \label{eqn:def_f2vw} \end{equation}
where
\[ \overline \Gamma_H(N) = \Gamma_H(N) / \Gamma(N) = \left\{ \mat{h^{-1}}{x}{0}{h} \in \SLZN \ \vert \ h \in H, x \in \Z/N\Z \right\}. \]
These forms do span~$\M_2(\Gamma_H(N))$, thanks to our assumption that~$p \nmid \# H$. Indeed,~$\overline \Gamma_H(N)$ has order~$N \# H$, and we have the following easy result:

\begin{lem}\label{lem:tr_surj}
Let~$M$ be a module over a ring~$R$, let~$G$ be a finite group of automorphisms of~$M$, and let
\[ M^G = \{ m \in M \ \vert \ \forall g \in G, \ g(m) = m \}. \]
If~$\# G$ is invertible in~$R$, then
\[ M^G = \left\{ \sum_{g \in G} g(m) \ \vert \ m \in M \right\}. \]
\end{lem}

\begin{proof}
The map
\[ \begin{array}{rcl} M & \longrightarrow & M \\ m & \longmapsto & \displaystyle \frac1{\# G} \sum_{g \in G} g(m) \end{array} \]
induces the identity on~$M^G$.
\end{proof}

As we explained, the Weierstrass model of~$E$ provides us with a normalisation of the differential on~$E$ and thus with a local trivialisation of~$\L$ at the corresponding point, so that the ``value'' of a modular form at this point is a well-defined quantity. Explicitly,~\eqref{eqn:fibre_XH} shows that given a primitive~$u \in (\Z/N\Z)^2$, the ``value'' of ~$f_{2,H}^{v,w}$ at the point of~$X_H(N)$ represented by~$\big(\W,H \cdot \beta(u)\big)$ is 
\begin{align*} f_{2,H}^{v,w}\big(\W,H \cdot \beta(u)\big) &= \sum_{\gamma \in \overline \Gamma_H(N)} f_1^{v \gamma}(\W,\beta_U) f_1^{w \gamma}(\W,\beta_U) \\
&= \sum_{\gamma \in \overline \Gamma_H(N)} f_1^{v \gamma U}(\W,\beta) f_1^{w \gamma U}(\W,\beta) \in \Z_q \end{align*}
where~$U$ is any element of~$\SL_2(\Z/N\Z)$ whose bottom row is~$u$.
We can thus compute this ``value'' in~$\Z_q/p^e$ thanks to proposition~\ref{pro:Kamal_log}, since we have determined the coordinates in~$\Z_q/p^e$ of the points of~$E[N]$ in the model~$\W$.

\bigskip

In order to obtain a basis of~$\M_2\big(\Gamma_H(N)\big)$, which has dimension
\[ d_2 = \dim \M_2\big(\Gamma_H(N)\big) = g + \nu_\infty, \]
we simply successively pick random pairs~$(v,w)$ of nonzero vectors of~$(\Z/N\Z)^2$, and form for each such pair the vector of ``values'' of the form~$f_{2,H}^{v,w}$ at all the points of the fibre, that is to say the vector of the~$f_{2,H}^{v,w}\big(\W,H \cdot \beta(u)\big)$ where~$u$ ranges over the primitive vectors of~$(\Z/N\Z)^2$ mod~$H$, until the reduction mod~$p$ of these vectors has rank~$d_2$. We then extract a basis, and thus obtain the matrix~$V$ (in the notation of definition~\ref{de:Mak_model}) for our~$p$-adic Makdisi model.

\begin{rk}
By sticking to the moduli interpretation of modular curves, we have thus managed to obtain a~$p$-adic Makdisi model for~$X_H(N)$ without requiring plane equations nor writing  down a single~$q$-expansion, merely by looking at the~$N$-torsion of just one elliptic curve~$E$ over~$\Q$. Besides, this method is straightforward to generalise to modular curves corresponding to any congruence subgroup. It could even be generalised to Shimura curves if an analogue was known for Makdisi's moduli-friendly forms in this context, but sadly this does not seem to be the case at the time of writing.
\end{rk}

\begin{rk}\label{rk:l1_l2_2}
As mentioned in remark~\ref{rk:l1_l2_1}, we may be tempted to construct our forms of weight 2 by taking products of two forms of the form~$l_1^{v,w} : (\W,\beta) \mapsto \lambda_{\beta(v),\beta(w)}$ rather than the~$f_1^v$. The matrix~$V$ has size~$O(g) \times O(g)$, so this would require us to evaluate~$O(g^2) O(\# \overline \Gamma_H(N))$ such forms~$l_1^{v,w}$, and therefore to perform that many divisions in~$\Z_q/p^e$. In the case where~$\Gamma_H(N) = \Gamma_0(N)$ with~$N$ prime, that is~$O(N^4)$ divisions in~$R$; whereas in the case where~$\Gamma_H(N) = \Gamma_1(N)$ with~$N$ prime, that is~$O(N^5)$ divisions. In contrast, if we work with the~$f_1^v$, we can precompute the~$N \times N$ matrix containing the~$f_1^v(\W,\beta)$ for all~$v \in (\Z/N\Z)^2$, which only requires~$O(N^2 \log N)$ divisions with the algorithm outlined in proposition~\ref{pro:Kamal_log}, after what no further divisions are required to fill in the matrix~$V$. Such a precomputation would not be of any help with the~$l_1^{v,w}$, since there are~$O(N^4)$ sets of the form~$\{v,w,-v-w\}$ with~$v$ and~$w$ in~$(\Z/N\Z)^2$, and therefore that many forms~$l_1^{v,w}$.
\end{rk}

\begin{rk}
In order to compute Galois representations by strategy~\ref{algo:Strategy_Hensel}, we need to be able to generate torsion points in the Jacobian over~$\F_q$, and for this, we must generate random points over~$\F_q$. In Makdisi's algorithms, a point on the Jacobian is represented by a subspace of a fixed Riemann-Roch space defined by vanishing conditions at an effective divisor on the curve. Bruin~\cite[Algorithm 3.7]{Bruin} presents a sophisticated method to generate uniformly distributed random points on the Jacobian in the framework, but as explained in~\cite[6.2.1]{Hensel}, we use a much cruder (and faster) approach, which in the case of modular curves amounts to considering subspaces of modular form spaces consisting of forms that vanish at the points of the modular curve represented by~$(\W,H\cdot\beta(u))$ for some randomly chosen~$\gamma \in \SLZN$. Since these points are in a rather special configuration, namely as they all lie on the same fibre of the projection to~$X(1)$, it may happen that the random points of~$J_H(N)(\F_q)$ obtained this way are so poorly distributed than they generate a subgroup with so little~$\ell$-torsion that it does not allow us to generate the representation space, so that the computation of the Galois representation stalls at stage~\ref{algo:Strategy_Hensel_randtors} of strategy~\ref{algo:Strategy_Hensel}. Fortunately, this seems rare in practice; in fact, in most of the cases that we have encountered, switching to another elliptic curve~$E$ suffices to solve this issue. Another workaround would consist in using not one but several elliptic curves~$E$, so as to allow ourselves to work with divisors supported on several fibres of the projection to~$X(1)$; and if this also fails, then we can fall back to Bruin's method.
\end{rk}

\subsection{The choice of~$p$}\label{sect:Choice_p}


Suppose we want to compute the mod~$\l$ representation~$\rho_{f,\l}$ attached to a newform~$f$ of weight~$k$, level~$N$, and nebentypus~$\eps_f$. As explained in section~\ref{sect:rho_in_XH}, this representation is found up to twist in the~$\ell$-torsion of~$J_H(N')$, where~$N'$ is defined by~\eqref{eqn:def_N'} and~$H$ is defined by~\eqref{eqn:def_H}.

As explained in the previous subsection, given a prime~$p \nmid 6\ell N \#H$, we can construct a~$p$-adic Makdisi model of~$J_H(N')$ by fixing a Weierstrass equation~$\W$ having good reduction at~$p$ and defining an elliptic curve~$E$ over~$\Q$ having~$j$-invariant neither 0 nor 1728 mod~$p$. This requires working in the unramified extension~$\Q_q = \Q_p(E[N])$, and in return allows us to compute explicitly with points of~$J_H(N')(\Z_q/p^e)$ for any~$e \in \N$ thanks to the methods presented in~\cite{Hensel}.

This leads to a method to compute~$\rho_{f,\l}$, provided that the points of~$J_H(N')[\ell]$ affording~$\rho_{f,\l}$ are defined over~$\Q_q$. This last requirement is equivalent to the degree~$a = [\Q_q:\Q_p]$ being a multiple of the order of~$\rho_{f,\l}(\Frob_p)$. This order can usually be determined explicitly, and in any case bounded, from the knowledge of the characteristic polynomial
\[ \chi_p(x) = x^2-a_p(f)x+ p^{k-1} \eps_f(p) \bmod \l \in \F_\l[x] \]
of~$\rho_{f,\l}(\Frob_p)$, as explained in proposition 6.1 of~\cite{Hensel}.

In summary, the degree~$a$ must satisfy two constraints, which both depend on~$p$: first, there must exist an elliptic curve as above having its~$N$-torsion defined over~$\Q_q$, which in particular imposes both
\begin{equation} q = p^a \equiv 1 \bmod N \label{eqn:Ep:Weil} \end{equation}
by the Weil pairing and
\begin{equation} q \geqslant (N-1)^2 \label{eqn:Ep:Hasse} \end{equation}  by the Hasse bound, and second, it must be a multiple of the order of~$\rho_{f,\l}(\Frob_p)$. Naturally, the smaller~$a$, the more efficient the computations will be (bearing in mind the remarks made in~\cite[6.4]{Hensel}), so it is a good idea to try many values of~$p$, and to select the one resulting in~$a$ being as small as possible. On the top of that,~$p$ must be such that~$\chi_p(x)$ is coprime mod~$\l$ with its cofactor in the L-factor~$L_p(x)$ of~$X_H(N')$ at~$p$ so that we can isolate the subspace of~$J_H(N')[\ell]$ affording~$\rho$, cf.~\eqref{eqn:chi_mult_1}.

A reasonable strategy is thus to determine in parallel~$\chi_p(x)$ and~$L_p(x)$ for all~$p$ not dividing~$6N' \#H$ up to some bound~$B$, and to retain a value of~$p$ leading to a degree~$a$ which is as small as possible. The value of~$B$ depends on how fast we can determine~$\chi_p(x)$, which involves evaluating~$a_p(f) \bmod \l$, and~$L_p(x)$, which by~\eqref{eqn:Lp_mod_res} involves computing the action of the Hecke operator~$T_p$ on~$\S_2\big(\Gamma_H(N')\big)$; in practice, we use~$B=100$ or~$1000$, cf. the examples in section~\ref{sect:Examples} below.

\begin{rk} 
Conditions~\eqref{eqn:Ep:Weil} and~\eqref{eqn:Ep:Hasse} are necessary, but not sufficient, for there to exist a suitable elliptic curve~$E$. For instance, there exists no elliptic curve over~$\F_{13}$ having~$j \not \in \{0,1728\}$ and full~$4$-torsion over~$\F_{13}$. As a result, extra care must be taken when trying small values of~$p$. Determining necessary and sufficient conditions on~$p$ and~$a$ in terms of~$N$ is an interesting problem, that could probably be solved by examining the Zeta function of the modular curve~$X(N)$; we have chosen not to go this way, and to simply try random Weierstrass equations until we find a curve having full~$N$-torsion over~$\F_q$ and~$j \not \in \{0,1728\}$, and to give up this value of~$p$ if no such curve is found after a certain number of attempts.
\end{rk}

\subsection{Finding a suitable elliptic curve and computing a basis of it~$N$-torsion}\label{sect:LiftEN}

Let~$N \in \N$ be an integer,~$p \nmid 6N$ a prime, and~$a \in \N$ a degree such that there exists an elliptic curve~$E$ as above, that is to say defined over~$\Q$, having good reduction at~$p$, having~$j$-invariant distinct from~$0$ and~$1728$ mod~$p$, and having all its~$N$-torsion defined over the unramified extension~$\Q_q$ of~$\Q_p$ of degree~$a$; in particular,~\eqref{eqn:Ep:Weil} and~\eqref{eqn:Ep:Hasse} must be satisfied. The purpose of this section is to explain how to find such a curve efficiently. The typical range that we have in mind is~$N \leqslant 1000$,~$p \leqslant 10^4$, and~$a \leqslant 100$.

Since~$p \nmid N$, the requirement that~$\Q_p(E[N]) \subseteq \Q_q$ is equivalent to~$\F_p(\bar E[N]) \subset \F_q$, where~$\bar E$ denotes the reduction of~$E$ mod~$p$; therefore, since~$p \nmid 6$, we will actually look for integers~$0 < A, B < p$ such that the short Weierstrass equation
\[ (\W) : y^2 = x^3 + Ax + B \]
viewed mod~$p$ defines an elliptic curve~$\bar E_{A,B}$ over~$\F_p$ (so that~$4A^3+27B^2 \not \equiv 0 \bmod p$), whose~$j$ invariant will automatically be distinct from~$0$ and~$1728$ mod~$p$ since~$A$ and~$B$ are nonzero mod~$p$. Our assumption is that there exists at least one such pair~$(A,B)$ such that the~$N$-torsion of this curve is defined over~$\F_q$, so that~$(\W)$ then defines an elliptic curve over~$\Q$ with the desired properties.

Our strategy simply consists in trying random pairs~$(A,B)$ until the condition 
\begin{equation} \F_p(\bar E_{A,B}[N]) \subseteq \F_q \label{eqn:cond_EAB} \end{equation}
is satisfied. Given such a random pair, we expect that~\eqref{eqn:cond_EAB} will most likely not be satisfied, so instead of directly testing~\eqref{eqn:cond_EAB} by computing the~$N$-division polynomial~$\psi_N(x)$ of~$\bar E_{A,B}$ which would be time-consuming, we begin by submitting~$\bar E_{A,B}$ to a battery of quick tests based on point-counting and aiming at weeding out most of the pairs~$(A,B)$ for which~\eqref{eqn:cond_EAB} does not hold. Once we find a pair~$(A,B)$ which passes these tests, we then submit it to extra tests to try to prove that~\eqref{eqn:cond_EAB} does hold, still while trying to avoid expensive computations such as the determination of~$\psi_N(x)$.

Our point is that since~$N$ is reasonably small, we can factor it as~$N = \prod_k l_k^{v_k}$ where the~$l_k \in \N$ are distinct primes, and then~\eqref{eqn:cond_EAB} is equivalent to
\[ \F_p(\bar E_{A,B}[l_k^{v_k}]) \subseteq \F_q \]
for all~$k$. If we let~$\Frob_p : x \mapsto x^p$ be the standard pro-generator of~$\Gal(\overline \F_p/\F_p)$, and if we define~$\Frob_q = \Frob_p^a : x \mapsto x^q$, then this can be rephrased by saying that~$\Frob_q$ must act trivially on~$\bar E_{A,B}[l_k^{v_k}]$ for all~$k$.

Since by assumption~$p$ is not too large, given a pair~$(A,B)$, we can quickly determine the quantity
\[ a_p = p+1-\# \bar E_{A,B}(\F_p) \in \Z. \]

The characteristic polynomial of~$\Frob_p$ acting on~$\bar E_{A,B}$ is then~$\chi(x) = x^2-a_p x +p \in \Z[x]$. Given~$a_p$, it is thus straightforward to compute its discriminant~$\Delta = a_p^2 - 4p \in \Z$, as well as the Newton sum
\[ \nu_a = \alpha^a + \beta^a \in \Z \]
where~$\alpha$ and~$\beta$ are the roots of~$\chi(x)$ in~$\overline \Q$ and~$a = [\F_q:\F_p]$ as above; this can even be done symbolically, without actually computing~$\alpha$ and~$\beta$. Naturally, if the cardinality
\[ \# \bar E_{A,B}(\F_q) = q+1-\nu_a \]
is not a multiple of~$N^2$, then the pair~$(A,B)$ can be rejected.

Define 
\[ M_1 = \prod_{l_k \nmid \Delta} l_k; \]
then the action of~$\Frob_p$ on~$\bar E_{A,B}[M_1]$ is semisimple. The characteristic polynomial of~$\Frob_q$ acting on~$\bar E_{A,B}[M_1]$ is~$(x-\alpha^a)(x-\beta^a) \equiv x^2-\nu_a x + 1 \in \Z/M_1\Z[x]$ since~$p^a \equiv 1 \bmod N$ by assumption, so~$\Frob_q$ acts trivially on~$E[M_1]$ iff.~$\nu_a \equiv 2 \bmod M_1$. We can therefore reject the pair~$(A,B)$ if this condition is not satisfied.

If now~$l_k$ is one of the prime factors of~$N$ dividing~$\Delta$, then~$\Frob_p \circlearrowright \bar E_{A,B}[l_k]$ is not semisimple, and therefore has a single eigenvalue~$c \in \F_{l_k}$ which satisfies~$2c \equiv a_p \bmod l_k$ as can been seen by considering the trace. If~$l_k=2$, then~$\Frob_p \circlearrowright \bar E_{A,B}[l_k]$ is necessarily unipotent, and therefore so is~$\Frob_q$; else, ~$\Frob_q \circlearrowright \bar E_{A,B}[l_k]$ is unipotent iff.~$a_p^a \equiv 2^a \bmod l_k$, therefore the pair~$(A,B)$ can be rejected if this condition is not satisfied for at least one of such~$l_k$.

These three simple tests eliminate most of the pairs~$(A,B)$. We now assume that~$(A,B)$ has passed these three tests, which means that the action of~$\Frob_q$ is trivial on~$\bar E_{A,B}[M_1]$ and unipotent on~$\bar E_{A,B}[l_k]$ for each~$l_k \mid \Delta$; it remains to determine whether~$\Frob_q$ really acts trivially on~$\bar E_{A,B}[l_k^{v_k}]$ for each~$k$. This is automatically the case for the~$l_k \nmid \Delta$ such that~$v_k=1$, as well as for the~$l_k \mid \Delta$ such that~$v_k=1$ and~$l_k \mid a$ since a unipotent mod~$l_k$ matrix of size~$2 \times 2$ which is also an~$a$-th power is then necessarily trivial; we therefore do not consider these primes anymore.

For each of the remaining primes, we then compute the division \linebreak polynomial~$\psi_{l_k^{v_k}}(x) \in \F_p[x]$ of~$\bar E_{A,B}$, and determine the degrees of its factors over~$\F_p$, which is faster than factoring it completely~\cite[3.4.3]{GTM138}. If these degrees do not all divide~$a$, then this polynomial does not split over~$\F_q$, so the action of~$\Frob_q$ on~$\bar E_{A,B}[l_k^{v_k}]/\pm1$ is nontrivial and the pair~$(A,B)$ can be rejected. Else, for each~$k$ such that~$l_k^{v_k}\neq 2$, we determine the roots of~$\psi_{l_k^{v_k}}(x)$ in~$\F_q$, and for each such root~$z$, we check whether~$z^3+Az+B$ is a square in~$\F_q$ by raising it to the~$\frac{q-1}2$ using fast exponentiation (if~$l_k^{v_k}=2$, then this will automatically be satisfied since~$\bar E_{A,B}[2]/\pm1 = \bar E_{A,B}[2]$). If this is the case, we have found a suitable pair~$(A,B)$.

\bigskip

Suppose now that we have found a suitable pair~$(A,B)$. In order to compute a~$p$-adic Makdisi model for~$X_H(N)$ to accuracy~$O(p^e)$, where~$ \in \N$ is a fixed parameter, we need to determine the coordinates in~$\Z_q/p^e$ of the~$N$-torsion points of the elliptic curve~$E_{A,B}$ over~$\Z_p$ defined by~$(\W)$. It is sufficient to determine the coordinates of two points~$P,Q$ forming a basis of~$E_{A,B}[N]$, since the coordinates of the other torsion points can then be obtained by applying the group law of~$E_{A,B}(\Z_q/p^e)$ as~$(\W)$ has good reduction at~$p$. We must also compute the matrix expressing how~$\Frob_p$ acts on~$E_{A,B}[N]$ with respect to this basis, so as to determine how~$\Frob_p$ permutes the points of the fibre of~$X_H(N) \longrightarrow X(1)$ corresponding to~$E$. Besides, later we will also need the value~$e_N(P,Q)$ of the Weil pairing of this basis, which is a primitive~$N$-th root of~$1$ in~$\Z_q/p^e$.

Again, we want to try to avoid the expensive computation of the~$N$-division polynomial of~$E_{A,B}$, so we proceed prime-by-prime. Factor as above~$N = \prod_k l_k^{v_k}$ where the~$l_k$ are distinct primes, define~$N_k = N / l_k^{v_k}$ or each~$k$, and let~$i_k \in \Z/N\Z$ be the idempotents corresponding to the Chinese remainder decomposition
\[ \Z/N\Z \simeq \prod_k \Z/l_k^{v_k}\Z, \]
that is to say
\[ i_k \bmod l_j^{v_j} = \left\{ \begin{array}{l} 1 \text{ if } j=k, \\ 0 \text{ if } j \neq k; \end{array} \right. \]
these~$i_k$ may be computed using B\'ezout relations between~$N_k$ and~$l_k^{v_k}$.

For each~$k$, we begin by computing the polynomials~$\psi_{l_k^{v_k}}(x)$ and~$\psi_{l_k^{v_k-1}}(x)$, \linebreak
where~$\psi_m(x) \in \Q[x]$ denotes the~$m$-th division polynomial of~$E_{A,B}$. We then pick two roots~$\bar x_{P_k}, \bar x_{Q_k}$ of~$\psi_{l_k^{v_k}}(x)$ in~$\F_q$, neither of which is not a root of~$\psi_{l_k^{v_k-1}}(x)$, and we set~$\bar P_k = (\bar x_{P_k}, \bar y_{P_k}), \ \bar Q_k = (\bar x_{Q_k}, \bar y_{Q_k})$, where~$\bar y_{P_k} \in \F_q$ is either square root of~$\bar x_{P_k}^3+A \bar x_{P_k} + B$, and similarly for~$\bar y_{Q_k}$. Then~$\bar P_k, \bar Q_k \in \bar E_{A,B}(\F_q)$ are two points of exact order~$l_k^{v_k}$; in particular, we can compute their Weil pairing
\[ \bar z_k = e_{l_k^{v_k}}(\bar P_k, \bar  Q_k) \in \F_q, \]
which is a primitive~$l_k^{v_k}$-root of~$1$ iff.~$\bar z_k^{l_k^{v_k-1}} \neq 1$. If this is not the case, then we start over with another choice of~$\bar x_{P_k}, \bar x_{Q_k}$; else we have obtained a basis of~$\bar E_{A,B}[l_k^{v_k}]$ over~$\F_q$. We now assume that this is the case.

We can then determine the matrix of~$\Frob_p$ acting on~$E_{A,B}[l_k^{v_k}]$ with respect to (the unique~$p$-adic lift of) this basis as
\[ \Phi_k = \mat{\log_{\bar z_k} e_{l_k^{v_k}}(\bar Q_k, \bar P_k^{\Frob_p})}{\log_{\bar z_k} e_{l_k^{v_k}}(\bar Q_k,\bar Q_k^{\Frob_p})}{-\log_{\bar z_k} e_{l_k^{v_k}}(\bar P_k,\bar P_k^{\Frob_p})}{-\log_{\bar z_k} e_{l_k^{v_k}}(\bar P_k,\bar Q_k^{\Frob_p})}, \]
where~$\log_{\bar z_k} : \mu_{l_k^{v_k}}(\F_q) \longrightarrow \Z/l_k^{v_k}\Z$ denotes the discrete logarithm in base~$\bar z_k$.

Next, we lift this basis~$(\bar P_k, \bar Q_k)$ from~$\F_q$ to~$\Z_q/p^e$ by first Hensel-lifting the~$x$-coordinates as roots of~$\psi_{l_k^{v_k}}(x)/ \psi_{l_k^{v_k-1}}(x)$, and then the~$y$-coordinates as square roots of~$x^3+Ax+B$; we thus obtain a basis~$(P_k, Q_k)$ of~$E_{A,B}[l_k^{v_k}]$ over~$\Z_q/p^e$. In principle, the value~$z_k \in \Z_q/p^e$ of its Weil pairing could also obtained by lifting~$\bar z_k$ as a root of~$x^{l_k^{v_k}}-1$, but we defer this for now since we will see that we can do better.

\begin{rk}
Some of the division polynomials~$\psi_{l_k^{v_k}}(x)$ may have been computed in~$\F_p[x]$ during the earlier phase when we searched for an appropriate pair~$(A,B)$; however they need to be re-computed, since we need their value in~$\Q[x]$ (as opposed to mod~$p$) here.
\end{rk}

It is then clear that
\[ P = \sum_k P_k, \ Q= \sum_k Q_k \]
form a basis of~$E_{A,B}[N]$, with respect to which the matrix~$\Phi \in \GLZN$ describing the action of~$\Frob_p$ can be obtained from the~$\Phi_k$ by Chinese remainders thanks to the idempotents~$i_k$. Besides, its Weil pairing
\begin{equation} z = e_N(P,Q) \in \mu_N(\Z_q/p^e) \label{eqn:zetaN} \end{equation}
may be determined off the~$\bar z_k$. Indeed, it is enough to determine~$\bar z = z \bmod p \in \F_q$, since we can then Hensel-lift this value as a root of~$x^N-1$. Furthermore, given an elliptic curve~$E$ and two integers~$M_1, M_2 \in \N$, the definition of the Weil pairing in terms of meromorphic functions on~$E$ with prescribed divisors shows that we have the identity
\[ e_{M_1M_2}(R,S) = e_{M_1}(R,S)^{M_2} \]
for all~$R,S \in E[M_1]$. Therefore, we find that
\begin{align*} \bar z &= \bar z^{\sum_k i_k} = \prod_k \bar z^{i_k} = \prod_k \bar z^{i_k^2} = \prod_k e_N(i_k \bar P, i_k \bar Q) \\ &= \prod_k e_N(\bar P_k,\bar Q_k) = \prod_k e_{l_k^{v_k}}(\bar P_k,\bar Q_k)^{N_k} = \prod_k \bar z_k^{N_k}. \end{align*}
The advantage of this approach is that the Weil pairing computations, which can be expensive, are only performed on the~$\bar E_{A,B}[l_k^{v_k}](\F_q)$ instead of~$E_{A,B}[N](\Z_q/p^e)$.

\subsection{Optimising the Makdisi model}\label{sect:pruning}

Makdisi's algorithms to compute in Jacobians rely on linear algebra involving matrices whose dimensions are determined by the parameters~$d_0$ and~$n_Z$ introduced in definition~\ref{de:Mak_model}. The smaller these parameters, the faster the computations; however these parameters must respectively satisfy the bounds~\eqref{eqn:d0bound} and~\eqref{eqn:nZbound} for Makdisi's algorithms to be valid.

The~$p$-adic Makdisi model of~$X_H(N)$ that we constructed in subsection~\ref{sect:Mak_model_strategy} above satisfies these bounds, and actually exceeds them, often quite significantly so. The purpose of this subsection is to show how to optimize it by tweaking it so that~\eqref{eqn:d0bound} and~\eqref{eqn:nZbound} are satisfied as sharply as possible, which in practice results in a major speedup of our computations.

Let us begin with~\eqref{eqn:d0bound}. Recall that we chose~$\L$ to be the line bundle whose sections are~$\M_2\big(\Gamma_H(N)\big)$, so that~$d_0 \overset{\text{def}}{=} \deg \L = 2g-2+\nu_\infty$. This ensured that the bound~\eqref{eqn:d0bound} on~$d_0$ is satisfied since we assume that the number~$\nu_\infty$ of cusps of~$X_H(N)$ is at least 3. We can thus make~\eqref{eqn:d0bound} an equality  as in~\cite[3.3]{algo}, that is to say by fixing three cusps of~$X_H(N)$ and by replacing~$\L$ with the sub-sheaf whose sections are the modular forms of weight 2 that vanish at all cusps except maybe these three.

In order to achieve this, we begin as in subsection~\ref{sect:Mak_model_strategy} by finding a basis of the space~$\M_2\big(\Gamma_H(N)\big)$ consisting of forms~$f_{2,H}^{v_i,w_i}$ as defined by~\eqref{eqn:def_f2vw} for various~$v_i,  w_i \in (\Z/N\Z)^2$. We then determine the value of these forms at each cusp except these three, and we deduce by linear algebra a basis of the subspace consisting of forms that vanish at all cusps except maybe these three.

In view of~\eqref{eqn:def_f2vw}, this requires determining the value of~$f_1^v$ at each cusp for each nonzero~$v \in (\Z/N\Z)^2$. Note that this value is well-defined by~\eqref{eqn:loc_param} applied to the case~$n=0$, and that it is actually enough to determine the value of~$f_1^v$ at the cusp~$\infty$ in terms of~$v$ since~$f_1^v \mact \gamma = f_1^{v \gamma}$ for all~$\gamma \in \SLZN$.

We are actually going to determine the whole~$q_N$-expansion of~$f_1^v$, since this will be useful in section~\ref{sect:Eval} below. By~\cite[corollary 3.13]{MakEis},~$f_1^v(\W,\beta)$ is proportional to the Weierstrass Zeta function of the elliptic curve defined by~$(\W)$ evaluated at the~$N$-torsion point~$\beta(v)$, and therefore to the modular form denoted by~$g_1^v$ in~\cite[4.8]{DS}. Combining the formulas found in~\cite{DS}, and in particular in section 4.8 thereof, we then find after some computations whose details we omit the following formula.
\begin{thm}\label{thm:qexp}
There exists a constant~$C$ depending only on~$N$ such that for all~$0 \leqslant c < N$ and~$d \in \Z/N\Z$ such that~$v=(c,d)$ is nonzero in~$(\Z/N\Z)^2$, we have
\[ f_1^v = C \sum_{n=0}^{+\infty} a_n q_N^n, \]
where
\[ a_0 = \left\{ \begin{array}{cc} \displaystyle \frac12 \frac{1+z^d}{1-z^d} & \text{ if } c=0, \\ \displaystyle \frac12 - \frac{c}N & \text{ if } c \neq 0, \end{array} \right. \]
and for all~$n \geqslant 1$,
\[ a_n = \sum_{\substack{a,b \in \Z \\ ab = n \\ a \equiv c \bmod N}} \sgn(b) z^{bd}, \]
where~$z=e_N(P,Q)$ is the primitive~$N$-th root of unity defined by~\eqref{eqn:zetaN} and which tells us which geometric component of~$X(N)$ we are working in.
\end{thm}

\begin{rk}
It should also be possible to derive these formulas by evaluating the moduli-friendly form~$f_1^v$ on the Tate curve. However, although this would be more in the spirit of this article, this seems to lead to more difficult computations.
\end{rk}

\begin{rk}
These formulas allow us to determine the coefficients~$a_n$ for~$n$ up to some bound~$B$ in quasi-linear time in~$B$. Since apart from cusp forms of weight~1, every modular form over a congruence subgroup of level~$N$ is expressible as a polynomial in the~$f_1^v$ by theorem~\ref{thm:MakEis}, we can thus compute the first~$B$ coefficients of the~$q$-expansion of any such form in quasi-linear time in~$B$ thanks to fast series arithmetic. This is quasi optimal, and faster than other methods such as modular symbols whose complexity is at least quadratic in~$B$; furthermore, by nature this approach is well-suited to  the Chinese remainder strategy which involves the computation of the desired result modulo several primes, and which is the key to many fast algorithms in computer algebra~\cite[5]{vzGathen}. However, the complexity of this approach with respect to the level~$N$ is probably terrible. Anyway, this is irrelevant for this article, since we will make very little use of~$q$-expansions, and since we will jut need a few terms when we do (cf. section~\ref{sect:Eval} below).
\end{rk}

We can thus optimise~$\L$ so that the bound~\eqref{eqn:d0bound} on~$d_0$ is sharp. Since~$d_0 = \deg \L$ drops, we can then reduce the number of points in the fiber of~$X_H(N) \longrightarrow X(1)$ at which we evaluate our forms. In fact, by~\eqref{eqn:nZbound}, we only need to retain~$10g+6$ of these points, since we now have~$d_0=2g+1$ exactly; however, the definition~\ref{de:Mak_model} of a Makdisi model also requires the set of these points to be globally invariant under~$\Frob_p$. Since we have determined how~$\Frob_p$ acts on the~$N$-torsion of the elliptic curve corresponding to this fibre, we can explicitly decompose this fibre into~$\Frob_p$-orbits; we then discard some of these orbits so that the bound~\eqref{eqn:nZbound} is satisfied and as sharp as possible, and we only evaluate our forms at the points in the remaining orbits to construct our~$p$-adic Makdisi model of~$X_H(N)$.

\begin{rk}
Note that unlike the bound~\eqref{eqn:d0bound} on~$d_0$, we cannot in general make the bound~\eqref{eqn:nZbound} on the number of points at which we evaluate our forms an equality, even though we can usually get pretty close since we have chosen the prime~$p$ so that the order~$a$ of~$\Frob_p$ is small.
\end{rk}

\section{Variant: Using the Hecke operator~$T_p$}\label{sect:Tp}

\subsection{Necessity of the use of~$T_p$}

Let~$f = q+\sum_{n\geq2} a_n(f) q^n \in \newf_k(N,\eps_f)$ be a newform, and suppose that we wish to compute the Galois representation~$\rho_{f,\l}$ attached to~$f$ mod a prime~$\l$ of degree~$1$ above~$\ell \in \N$. Let~$N'$ and~$H$ be as in~\eqref{eqn:def_N'} and~\eqref{eqn:def_H}, so that~$\rho_{f,\l}$ occurs in~$J_H(N')[\ell]$.

A limitation of the~$p$-adic strategy~\ref{algo:Strategy_Hensel} is that it assumes the existence of a good prime~$p$ satisfying condition~\eqref{eqn:chi_mult_1}, in other words such that the characteristic polynomial
\[ \chi_p(x) = x^2-a_p(x)+p^{k-1} \eps_f(p) \bmod \l \in \F_\l[x] \]
of~$\rho_{f,\l}(\Frob_p)$ is coprime mod~$\ell$ with the local L factor
\[ L_p(x) = \det\big(x-\Frob_p \vert_{J_H(N')}\big) \in \Z[x]. \]
In some rare cases, it may happen that condition~\eqref{eqn:chi_mult_1} is not satisfied by any prime~$p$, so that our~$p$-adic method to compute~$\rho_{f,\l}$, as presented so far, does not apply.

\begin{ex}\label{ex:Frobp_cannot_cut}
This happens when~$f$ is the newform of~\cite{LMFDB} label \mfref{5.6.a.a} and~$\l = 13$. This phenomenon is related to~$f$ being supersingular at~$\l$ and having trivial nebentypus. 
\end{ex}

However, in the immense majority of cases, including that mentioned in example~\ref{ex:Frobp_cannot_cut} above, multiplicity one statements such as~\cite[theorem 3.5]{RS} show that this can be remedied by isolating the subspace~$T_{f,\l}$ of dimension~$2$ of~$J_H(N')[\ell]$ affording~$\rho_{f,\l}$ as the subspace where the Hecke algebra acts with the eigenvalue system of~$f \bmod \l$ instead of in terms of the action of the Frobenius at a good prime~$p$. In other words, the representation space~$T_{f,\l}$ can be carved out as
\begin{equation} T_{f,\l} = \bigcap_{n \in \N} \ker(T_n - a_n(f) \bmod \l \vert_{J_H(N')[\ell]}). \label{eqn:cut_by_Tn} \end{equation}

What happens in example~\ref{ex:Frobp_cannot_cut} is that the \emph{generalised} kernels
\[ \bigcap_{n \in \N} \ker\big((T_n - a_n(f) \bmod \l)^\infty \vert_{ J_H(N')[\ell]}\big) \]
yield an~$\Fl$-subspace of dimension~4 which is a non-split extension of one copy of~$\rho_{f,\l}$ by another. This explains why an attempt based on the characteristic polynomials~$\chi_p(x)$ alone cannot succeed in this case.

In order to remedy this situation,~\eqref{eqn:cut_by_Tn} thus suggests we implement the action of the Hecke algebra on Makdisi models of modular curves. Although this is theoretically possible since~\cite{Bruin} shows that pull-backs and push-forwards are computable in Makdisi models, this approach seems complicated, and we have chosen not to follow it. Instead, we hope that there exists a prime~$p$ such that the representation space can be carved out simply as
\begin{equation} T_{f,\l} = \ker(T_p - a_p(f) \bmod \l)_{\vert J_H(N')[\ell]}, \label{eqn:cut_by_Tp} \end{equation}
which is tantamount to having~$a_p(f') \bmod \l' = a_p(f) \bmod \l$ for a \emph{unique} pair~$(f',\l')$ with~$f'$ an eigenform of weight~$2$ and level~$\Gamma_H(N')$ and~$\l'$ a prime of~$K_{f'}$ above~$\ell$. We have not encountered any case where no such~$p$ exists, so this seems to be a reasonable assumption.

\subsection{Implementing~$T_p$ on~$p$-adic Makdisi models}

The Eichler-Shimura relation \cite[8.7.2]{DS} states that~$T_p = \Frob_p + p \langle p \rangle_* \Frob_p^{-1}$ on~$J_H(N')(\overline \F_p)$. Therefore, if we construct a~$p$-adic Makdisi model of~$J_H(N')$ with the same prime~$p$ as in~\eqref{eqn:cut_by_Tp} so we may apply~$\Frob_p$, and which contains the extra data required to apply~$\langle p \rangle_*$, then we get an implementation of~$T_p$ on~$J_H(N')(\overline \F_p)$, and we may thus alter strategy~\ref{algo:Strategy_Hensel} to carve out the representation space with~$T_p$ instead of~$\Frob_p$.

In order to construct such a~$p$-adic Makdisi model, we must first find a good prime~$p$ such that~\eqref{eqn:cut_by_Tp} holds. We follow a search procedure similar to that described in subsection~\ref{sect:Choice_p}, where we substitute~\eqref{eqn:cut_by_Tp} for the condition that~$\chi_p$ be coprime with its cofactor in~$L_p$. In order to test~\eqref{eqn:cut_by_Tp} for a given prime~$p$, we could compute the matrix of~$T_p$ (with respect to any~$\Z$-basis) acting on cuspidal modular symbols of level~$\Gamma_H(N')$, reduce it mod~$\ell$, and see if its~$a_p(f) \bmod \l$-eigenspace has dimension~$2$ and not more. Alternatively, in view of the pairing between cusp forms and cuspidal modular symbols, we can also compute the matrix of~$T_p$ acting on the space~$\S_2\big(\Gamma_H(N')\big)$ of cusp forms with respect to a basis of~$\S_2\big(\Gamma_H(N')\big)$ such that this matrix has integers entries, which can be done with \cite{gp}, reduce it mod~$\ell$, and see if its~$a_p(f) \bmod \l$-eigenspace has dimension~$2$ and not more. Indeed, although it is conceivable that the basis of~$\S_2\big(\Gamma_H(N')\big)$ chosen by~\cite{gp} does not remain a basis after reduction mod~$\ell$, or that the pairing between modular forms and modular symbols degenerates mod~$\ell$, these phenomena can only increase the dimension of the~$a_p(f) \bmod \l$-eigenspace, so we do get a sufficient criterion for~\eqref{eqn:cut_by_Tp} to be satisfied. The advantage of working with cusp forms rather than modular symbols is that we will need anyway to determine the matrix of~$T_p$ acting on cusp forms so as to compute~$L_p(x)$ by~\eqref{eqn:Lp_mod_res}.

Next, our~$p$-adic Makdisi model must satisfy extra requirements so as to be able to apply the diamond operator~$\langle p \rangle_*$, as explained in subsection~\ref{subs:JAuts}. We can easily determine the permutation induced by~$\langle p \rangle$ on the fibre of~$X_H(N') \longrightarrow X(1)$ thanks to~\eqref{eqn:Diam_on_fibre}. However, we must modify the optimisation process described in subsection~\ref{sect:pruning}, for two reasons. First, the line bundle~$\L$ must be invariant by~$\langle p \rangle$, so instead of taking the bundle whose sections are the forms of weight~$2$ that vanish at all but three cusps, we must take the bundle whose sections are the forms of weight~$2$ that vanish at all cusps outside a set~$S$ containing at least three cusps and which is stable not only by~$\GQ$, but also by~$\langle p \rangle$. Therefore, we may need to take a set~$S$ with slightly more than three elements, which makes the~$d_0 = \deg \L$ larger and thus takes us away from attaining sharpness in the bound~\eqref{eqn:d0bound}. Second, while we can still drop some of the points of the fibre of~$X_H(N') \longrightarrow X(1)$, we must ensure not only that the bound~\eqref{eqn:nZbound} is satisfied in spite of the accretion of~$d_0$, but also that the remaining points are globally invariant not only under~$\Frob_p$, but also under~$\langle p \rangle$. As a result, it may well be that we are unable to optimise our~$p$-adic Makdisi model as well as before, in which case all our computations in~$J_H(N')$ will be slower. For this reason, it is preferable to carve out~$T_{f,\l}$ with~$\Frob_p$ as before when possible, and to reserve this new approach to cases such as example~\ref{ex:Frobp_cannot_cut}. 

\begin{rk} In principle, it would be possible to construct \emph{two}~$p$-adic Makdisi models of~$X_H(N')$, namely a ``large one'' with the extra data for~$\langle p \rangle_*$ and a ``small'', better-optimised one without. We would then use the large model to generate points of~$T_{f,\l}$, convert them to points of the small model, and proceed with the small model for the rest of the computation. However, we have not yet implemented this conversion process.
\end{rk}

\subsection{Isolating~$T_{f,\l}$ by~$T_p$}

We now describe in detail our new approach to compute~$\rho_{f,\l}$, given a good prime~$p$ satisfying~\eqref{eqn:cut_by_Tp}. Although this approach remains valid even in the case where~$\chi_p(x) \, \Vert \, L_p(x)$, meaning that~$\chi_p(x)$ is coprime with its cofactor~$L_p(x)/\chi_p(x)$ so that we could carve out~$T_{f,\l}$ by using~$\Frob_p$, we are chiefly interested in the case where~$\chi_p(x) \nparallel \, L_p(x)$.

As explained in section~\ref{sect:Choice_p}, we can determine from~$\chi_p(x)$ an integer~$a \in \N$ such that the points of~$T_{f,\l}$ are defined over the extension~$\F_q$ of~$\F_p$ of degree~$a$ and that~$\F_q$ contains the~$\ell$-th roots of unity. We then construct a~$p$-adic Makdisi model of~$J_H(N')$ over~$\Q_q$ with high-enough~$p$-adic accuracy and which contains the extra data needed to apply~$\langle p \rangle_*$, as explained above. We then reduce this model mod~$p$ so as to obtain a Makdisi model of~$J_H(N)(\F_q)$, while retaining the high~$p$-adic accuracy model for later use. Writing~$J=J_H(N')(\F_q)$ for brevity from now on, we can then apply~$T_p$ on~$J$ as
\begin{equation} T_p = \Frob_p + p \langle p \rangle_* \Frob_p^{a-1}. \label{eqn:Tp_JFq} \end{equation}

\begin{rk}
We will actually need to apply~$T_p$ to points of~$J[\ell]$ only. Therefore, we can replace~\eqref{eqn:Tp_JFq} with
\[ T_p = \Frob_p + m \langle p \rangle_* \Frob_p^{a-1} \]
where~$m$ is any integer satisfying~$m \equiv p \bmod \ell$. In particular, working with large values of~$p$ does not slow down the application of~$T_p$, so there is no harm in choosing a large value of~$p$ so as to make~$a$ smaller, thus making the calculations faster.
\end{rk}

Let~$\upsilon_p(x) = \gcd(L_p(x),\chi_p(x)^\infty) \in \Fl[x]$ be the largest mod~$\ell$ factor of~$L_p(x)$ whose irreducible factors all divide~$\chi_p(x)$; by construction,~$\chi_p(x) \mid \upsilon_p(x) \, \Vert \, L_p(x)$,
 so we know how to generate random points of the subspace
\[ U = \ker \upsilon_p(\Frob_p) \subseteq J[\ell] \]
by using the action of~$\Frob_p$. Observe that~$U \supseteq T_{f,\l}$ is an~$\Fl$-space of dimension 
\[ d = \deg \upsilon_p(x) \]
since~$\upsilon_p(x) \, \Vert \, L_p(x)$.

As~$\F_q$ contains the~$\ell$-th roots of unity by assumption, algorithm 6 of \cite{Hensel} allows us to evaluate the Frey-R\"uck pairing
\[ [ \cdot, \cdot] : J[\ell] \times J/\ell J \longrightarrow \F_q^\times / \F_q^{\times \ell} \overset{\sim}{\longrightarrow} \Fl, \]
where the rightmost arrow consists in
\[ \begin{array}{ccc} \F_q^\times / \F_q^{\times \ell} & \overset{\sim}{\longrightarrow} & \mu_\ell(\F_q) \\ x & \longmapsto & x^{(q-1)/\ell} \end{array} \]
followed by the discrete logarithm with respect to some fixed primitive~$\ell$-th root of unity. This paring is perfect, so we will use it to detect~$\Fl$-linear dependency in~$J[\ell]$, similarly to algorithm 13 of \cite{Hensel}.

Consider now the following algorithm.

\begin{figure}[H]
\fbox{\parbox{\textwidth}{
\begin{enumerate}
\item Generate random points~$y_1,\cdots,y_d$ of~$J[\ell]$. Initialise a matrix~$M$ of size~$0 \times 0$ over~$\Fl$, and a matrix~$P$ of size~$d \times 0$ over~$\Fl$. Finally, set~$r \la 0$.
\item \label{algo:CarveTp_NewPt} Set~$d_+ \la 0$, and generate a random point~$x$ of~$U$.
\item \label{algo:CarveTp_Tpclosure} Set~$P' \la$ the matrix obtained by sticking the column containing~$[x,y_1],\cdots,[x,y_d]$ to the right of~$P$.
\item \label{algo:CarveTp_Indep} If~$P'$ has rank~$r+1$, then set~$d_+ \la d_+ + 1$,~$r \la r+1$,~$x_{r} \la x$, and~$P \la P'$. Then set~$x \la T_p x$, and go to step~\ref{algo:CarveTp_Tpclosure}.
\item \label{algo:CarveTp_PseudoRel} Else, let~$\lambda_1, \cdots, \lambda_{r+1} \in \Fl$ be the coefficients of a non-trivial linear dependency relation between the columns of~$P'$.
\item \label{algo:CarveTp_FakeRel} If~$\sum_{i=1}^r \lambda_i x_i + \lambda_{r+1} x$ does not evaluate to~$0 \in J$, then first set~$d_+ \la d_+ + 1$,~$r \la r+1$,~$x_{r} \la x$, and~$P \la P'$. Next, generate random points~$y \in J$ until the row~$[x_1,y],\cdots,[x_r,y]$ is linearly independent from the rows of~$P'$, find an~$i \leqslant d$ such that the~$i$-th row of~$P$ is in the~$\Fl$-span of the other rows of~$P$, and replace this~$i$-th row with~$[x_1,y],\cdots,[x_r,y]$ and~$y_i$ with~$y$. Finally, set~$x \la T_p x$, and go to step~\ref{algo:CarveTp_Tpclosure}.
\item If~$\sum_{i=1}^r \lambda_i x_i + \lambda_{r+1} x =0 \in J$ and $d_+ = 0$, go to step~\ref{algo:CarveTp_NewPt}.
\item \label{algo:CarveTp_TrueRel} If~$\sum_{i=1}^r \lambda_i x_i + \lambda_{r+1} x =0 \in J$ but $d_{+} \neq 0$, set~$\lambda_i \la - \lambda_i / \lambda_{r+1}$ for all~$i \leqslant r$, so that~$x = \sum_{i=1}^r \lambda_i x_i$. Replace~$M$ with
\[ \left(
\raisebox{0.5\depth}{
 \xymatrix @R=0.7pc @C=1pc{
 &&& 0\ar@{..}[rrr] \ar@{..}[dd] &&& 0 \ar@{..}[dd] & \lambda_1 \ar@{-}[ddddddd] \\
& M && \\
&&& 0 \ar@{..}[rrr] &&& 0 \\
0 \ar@{..}[rr] \ar@{..}[dddd] && 0 \ar@{..}[dddd]  & 0 \ar@{..}[dddrrr] \ar@{..}[rrr] &&& 0 \ar@{..}[ddd] \\
&&&1 \ar@{-}[dddrrr] &&& \\
&&&0 \ar@{..}[dd] \ar@{..}[ddrr]&&&\\
&&&&&& 0\\
0  \ar@{..}[rr]&& 0 & 0\ar@{..}[rr]&&0& 1 & \lambda_r \\
}} \right) \]
where the block with the sub-diagonal of ones has~$d_+$ rows, and then determine~$K = \ker(M-a_p(f) \bmod \l)$. If~$\dim K = 2$, then let~$(\kappa_{1,j}, \cdots, \kappa_{r,j})$~$(j \in \{1,2\})$ be a pair of vectors of~$\F_\ell^r$ forming a basis of~$K$, let~$b_j = \sum_{i=1}^r \kappa_{i,j} x_i$, and return the pair~$(b_1,b_2)$ of points of~$J[\ell]$. Else, go to step~\ref{algo:CarveTp_NewPt}.
\end{enumerate}
}}
\caption{Carving out~$T_{f,\l}$ by~$T_p$}
\label{algo:CarveTp}
\end{figure}

Every time we enter step~\ref{algo:CarveTp_NewPt}, we have that the~$r < d$ points~$x_1,\cdots,x_r \in U$ are linearly independent over~$\Fl$ and span a subspace~$X \subseteq U$ which is stable under~$T_p$, that~$M$ is the matrix of~$T_p$ on~$X$ with respect to~$x_1,\cdots,x_r$, and that the~$d \times r$ matrix~$P$ of pairings~$[x_j,y_i]$ has full rank~$r$. Then at step~\ref{algo:CarveTp_Tpclosure}, the~$d\times(r+1)$ matrix~$P'$ has rank either~$r+1$ or~$r$. If~$P'$ has rank~$r+1$, then~$x$ is linearly independent from the~$x_j$, so at step~\ref{algo:CarveTp_Indep} we append~$x$ to the~$x_j$, increase~$r$, and start over with~$T_p x$ instead of~$x$. If~$P'$ has rank~$r$, then at step~\ref{algo:CarveTp_PseudoRel} the~$\lambda_i$ are unique up to scaling, and either~$x$ is genuinely linearly dependent on the~$x_j$, or the linear forms~$[\cdot,y_i]$ do not separate the points of~$X$; we determine which alternative we are in by checking whether~$\sum_{i=1}^r \lambda_i x_i + \lambda_r x$ is zero or not. If it is not, then at step~\ref{algo:CarveTp_FakeRel}, we append~$x$ to the~$x_j$ thus increasing~$r$, then we modify one of the~$y_i$ so that the linear forms~$[\cdot,y_i]$ separate the points in the span~$X$ of the~$x_j$, and finally we start over with~$T_p x$ instead of~$x$. If~$\sum_{i=1}^r \lambda_i x_i + \lambda_r x = 0$, then~$x$ is linearly dependent on~$x_1,\cdots,x_r$ since the latter are linearly independent; in particular,~$\lambda_{r+1} \neq 0$. At this stage,~$x_1,\cdots,x_r$ span a~$T_p$-stable~$\Fl$-subspace~$X$ of~$U$, whose dimension~$r$ has increased by~$d_+$ since the last time we entered step~\ref{algo:CarveTp_NewPt}. So if~$d_+=0$, then the point~$x$ generated at step~\ref{algo:CarveTp_NewPt} was already in~$X$ at that time, so we simply try again with another random~$x \in U$. Else, at step~\ref{algo:CarveTp_TrueRel} we update~$M$, and see whether~$X$ contains the~$2$-dimensional subspace~$T_{f,\l} = \ker (T_p-a_p(f) \bmod \l)$. If it does, then we return a basis of~$T_{f,\l}$; else we go back to step~\ref{algo:CarveTp_NewPt} so as to enlarge the~$\Fl[T_p]$-module~$X$ by throwing in a new random point~$x \in U$.

We thus obtain a pair of points~$(b_1,b_2)$ of~$J[\ell]$ forming an~$\Fl$-basis of the representation space~$T_{f,\l}$, so we may proceed with the calculation of~$\rho_{f,\l}$ by our usual strategy~\ref{algo:Strategy_Hensel} from step~\ref{algo:Strategy_Hensel:Lift} on.

\begin{rk}
Some later steps of strategy~\ref{algo:Strategy_Hensel} actually require the matrix of~$\Frob_p$ and of~$\langle p \rangle_*$ on~$T_{f,\l}$ with respect to our basis~$(b_1,b_2)$. We can easily obtain the pairings~$[b_j,y_i]$ by taking linear combinations of the columns of the matrix~$P$ in step~\ref{algo:CarveTp_TrueRel} of algorithm~\ref{algo:CarveTp} in the same way as for the points~$x_j$, and then compare these pairings with the~$[\Frob_p b_j, y_i]$ so as to deduce the matrix of~$\Frob_p$; as for the matrix of~$\langle p \rangle_*$, it is simply the scalar matrix~$\eps_f(p) \bmod \l$.
\end{rk}

An explicit example of use of this method is presented on page~\pageref{ex:Tp} below.

\section{Construction of evaluation maps}\label{sect:Eval}

In order to be able to compute Galois representations following strategy~\ref{algo:Strategy_Hensel}, we still need to construct one or more rational maps \[ \alpha : J_H(N) \dashrightarrow \A^1\] defined over~$\Q$.

As in~\cite[2.2.3]{Hensel}, we begin by constructing a rational map \[ \alpha_{\PP} : J_H(N) \dashrightarrow \PP(V_2), \] where~$V_2$ denotes the space of global sections of~$\L^{\otimes 2}$. As explained in the same reference, this requires picking two linearly-inequivalent effective divisors~$E_1 \not \sim E_2$ on~$X_H(N)$ of degree~$d_0-g$, such that we can compute the subspace of~$V_2$ formed by sections that vanish at~$E_i$ for each~$i \in \{1,2\}$; besides, these divisors must be defined over~$\Q$ so as to ensure that~$\alpha_{\PP}$ is defined over~$\Q$.

For these reasons, we choose~$E_1$ and~$E_2$ so that they are supported by cusps, possibly with multiplicities. The rationality  condition is then easy to satisfy since the modular curve~$X_H(N)$ tends to have plenty of rational cusps as mentioned in remark~\ref{rk:cusps_XH_rat}; besides, as explained in remark~\ref{rk:whybounds},~$V_2$ is spanned by products to two forms~$f_{2,H}^{v_i,w_i}$, whose~$q$-expansions can be determined by~\eqref{eqn:def_f2vw} and theorem~\ref{thm:qexp}, so we can determine the subspaces of~$V_2$ corresponding to these divisors by linear algebra, even if these divisors have multiplicities.

\bigskip

In order to get an~$\A^1$-valued Galois-equivariant map, it remains to construct a rational map
\[ \PP(V_2) \dashrightarrow \A^1 \]
defined over~$\Q$. For this, we offer two strategies.

\subsection{Strategy 1: Using~$q$-expansions}

The first strategy involves~$q$-expansions; namely, as in~\cite[3.6]{algo}, we construct this map as
\begin{equation} \begin{array}{ccc} \PP(V_2) & \dashrightarrow & \A^1 \\ f & \longmapsto & \displaystyle \frac{a_{n_1}(f \mact M_1)}{a_{n_2}(f \mact M_2)}, \end{array} \label{eqn:eval_qexps} \end{equation}
where~$n_1, n_2$ are nonnegative integers and~$M_1, M_2 \in \SLZN$ yield rational~$q$-expansions in the sense of definition~\ref{de:Qqexp}. Indeed, recall that the elements of~$V_2$ are modular forms (of weight~$4$) by our choice of~$\L$. 

\begin{rk}
In subsection~\ref{sect:pruning}, we actually redefined~$\L$ so that its global sections are the forms of weight~$2$ that vanish at all cusps except three. Therefore, if the cusp~$M_1 \cdot \infty$ is not one of these three, then~$n_1$ should be at least~$2$; similarly for~$M_2$ and~$n_2$.
\end{rk}

As noted in remark~\ref{rk:always_1_Qcusp}, there always exists at least one matrix which yields rational~$q$-expansions, namely~$M = \smat{0}{1}{-1}{0}$. Therefore, this construction applies to every level, e.g. by taking~$M_1=M_2=M$ and~$n_1 \neq n_2$ (lest we get a constant map). In fact, there are always infinitely many possible choices for the parameters~$n_1$ and~$n_2$, and usually many choices for~$M_1$ and~$M_2$ as well. For instance, we can enumerate the cusps of~$X_H(N)$ for which there exists a matrix which yield rational~$q$-expansions, and try all the pairs~$(M_1,M_2)$ of (not necessarily distinct) matrices in this list and all integers~$n_1,n_2$ up to some bound~$B$, e.g.~$B = 5$.

We thus obtain several evaluation maps~$\alpha$, some of which may not be injective on the~$\Fl$-subspace of~$J_H(N)[\ell]$ which affords our Galois representation, and therefore not useful for our purpose; however, in practice, if the bound~$B$ is large enough, we always get many injective versions of~$\alpha$, and thus many versions of the polynomial~$F(x)$ which describes the Galois representation (cf. strategy~\ref{algo:Strategy_Hensel}). We then simply keep the ``prettiest'' version, for instance that having the smallest arithmetic height.

\begin{rk}
Practically, we should record the~$q$-expansions of forms~$f_i \in V_2$ forming a basis of~$V_2$ during the creation of the~$p$-adic Makdisi model of~$X_H(N)$: this makes evaluating~\eqref{eqn:eval_qexps} at~$f \in V_2$ easy, since we merely have to identify by linear algebra~$f$ as a linear combination of the~$f_i$ from their values at some points of the fibre of~$X_H(N) \rightarrow X(1)$.
\end{rk}

\subsection{Strategy 2: Using forms defined over~$\Q$}

Another strategy consists in fixing a basis~$(f_i)_i$ of the space~$V_2$, and in considering the map
\[ \begin{array}{ccc} \PP(V_2) & \dashrightarrow & \A^1 \\ \displaystyle \sum_i \lambda_i f_i & \longmapsto & \displaystyle \frac{\lambda_{i_1}}{\lambda_{i_2}}, \end{array} \]
where~$i_1,i_2 \leqslant \dim V_2$ are fixed integers. Again, there many choices of pairs~$(i_1,i_2)$, so that we get many versions of~$\alpha$ and of~$F(x)$.

This is in fact the adaptation of the approach that we used in~\cite[2.2.3]{Hensel}, which applies to any curve, not just modular curves. Its advantage is thus that it completely dispenses us of~$q$-expansion computations; however, the basis~$(f_i)_i$ of~$V_2$ must be made up of forms which are defined over~$\Q$ for the resulting map to be Galois-equivariant.

For this, it is enough to construct a basis of the section space of~$\L$ formed of forms which are defined over~$\Q$, since we can then generate~$V_2$ by taking products of two such forms (cf. remark~\ref{rk:whybounds}). The sections~$f_{2,H}^{v,w}$ of~$\L$ introduced in~\eqref{eqn:def_f2vw} are, in general, only defined over the cyclotomic field~$\Q(\mu_N)$. However, given a form~$f$ defined over~$\Q(\mu_N)$,~\eqref{eqn:Gal_cyclo_qexp} applied to the quotient~$f/f_0$, where~$f_0$ is a form defined over~$\Q$ and of the same weight as~$f$, combined with the fact that~$M = \smat{0}{1}{-1}{0}$ yields rational coefficients, shows that
\[ f^{\sigma_x} = f \mact M \smat{1}{0}{0}{x} M^{-1} = f \mact \smat{x}{0}{0}{1} \]
for all~$x \in \ZNX$, where~$\sigma_x \in \Gal\big(\Q(\mu_N)/\Q\big)$ is as in~\eqref{eqn:Gal_cyclo}. Therefore, for all nonzero~$v \in (\Z/N\Z)^2$ and~$w \in (\Z/N\Z)^2$ and for all~$y \in \Z/N\Z$, the section
\begin{equation} \sum_{x \in \ZNX} \zeta_N^{xy} \ f_{2,H}^{v,w} \mact \smat{x}{0}{0}{1} = \sum_{x \in \ZNX} \zeta_N^{xy} \sum_{\gamma \in \overline \Gamma_H(N)} f_1^{v \gamma \smat{x}{0}{0}{1}} f_1^{w \gamma \smat{x}{0}{0}{1}} \label{eqn:def_f2vwQy} \end{equation}
of~$\L$ is defined over~$\Q$
. Furthermore, under the assumption that~$p \nmid \phi(N)$ (a mild strengthening of the assumption~$p \nmid \#H$ made earlier in subsection~\ref{sect:Mak_model_strategy}), lemma~\ref{lem:tr_surj} shows that these sections generate the section space of~$\L$ over~$\Z_p/p^e$ for any~$e \in \N$.

\bigskip

Unfortunately, a bit of experimenting with our implementation has revealed that this approach results in polynomials~$F(x)$ whose arithmetic height is tremendously larger than those obtained with the approach based on~$q$-expansions, and which therefore require ridiculously high $p$-adic accuracy in order to be identified as elements of~$\Q[x]$, cf. remarks~\ref{rk:Qf_bad_1} and~\ref{rk:Qf_bad_2} below. This expresses the fact that products of two forms of the form~\eqref{eqn:def_f2vwQy} do not form ``nice''~$\Q$-bases of~$\M_4\big(\Gamma_H(N)\big)$, and that the linear algebra used in subsection~\ref{sect:pruning} to find forms which vanish at all but three cusps makes things even worse. For this reason, we only use Strategy 1 in practice.

\section{Comparison with the complex-analytic method and results}\label{sect:Examples}

\subsection{Examples of computations}

We conclude be giving some examples so as to demonstrate the performance of our implementation of the method presented in this article. In these examples, the newforms are specified by their~\cite{LMFDB} label.

When analysing these examples, the reader should bear in my mind that the difficulty of the computation of a mod~$\ell$ representation is governed by two essential parameters: the genus of the modular curve used in the computation of course, but also the number~$\ell$ itself, since the computation must process~$\#(\Fl^2 \setminus \{0\}) = \ell^2-1$ torsion points, and outputs a polynomial~$F(x) \in \Q[x]$ of degree~$\ell^2-1$, whose arithmetic height is likely to grow with~$\ell$, thus requiring more~$p$-adic accuracy.

\subsubsection{``Small'' examples}

We begin with three ``small'' examples. Unless explicitly stated otherwise, the times we give are the ones obtained by executing these examples on the author's laptop, which has 4 hyperthreaded cores. As our implementation makes heavy use of the fact that certain steps of the computation are easily parallelisable, we express the computation times as ``X seconds of CPU time, and Y seconds of real time". This does not mean that the computation took X+Y seconds, but that the computation took Y seconds, during which the cumulated CPU time (taking parallelisation into account) was X seconds.

\subsubsection*{A form of weight 2 and level 16}

The form
\[ f = \mfref{16.2.e.a} = q + ( -1 - i ) q^{2} + ( -1 + i ) q^{3} + O(q^{4}) \]
is up to Galois-conjugacy the only newform of weight~$k=2$ and level~$\Gamma_1(16)$. Since its coefficient field~$K_f = \Q(i)$ is an extension of~$\Q$ of degree~$2$, the modular curve~$X_1(16)$ is of genus~$g=2$. Since~$f$ is of weight~2, the Jacobian~$J_1(16)$ of~$X_1(16)$ contains the mod~$\l$ representation~$\rho_{f,\l}$ attached to~$\l$ for any prime~$\l$ of~$K_f$. For this example, let us take~$\l=(5,i-2)$, one of the two primes of~$K_f$ above~$5$.

As explained in subsection~\ref{sect:Mak_model_strategy}, we must begin by choosing a prime~$p$ to work with. After trying all primes up to~$100$, which requires computing~$a_p(f)$ for~$p \leqslant 50$, we decide to take~$p=23$, because~$\rho_{f,\l}(\Frob_{23})$ has order only~4. This search takes about 640ms of CPU time, but only 110ms of real time, thanks to parallelisation.

Next, we construct a~$23$-adic Makdisi model of~$X_1(16)$ with residue degree~$a=4$ and accuracy~$O(23^{e})$, where we have chosen $e=7$ so as to be able to identify rationals of height at most $4 \times 10^4$. This construction involves spotting the elliptic curve
\[ y^2 = x^3+3x+3 \]
which has all its~$16$-torsion defined over the degree-4 unramified extension of~$\Q_{23}$. All this takes only 120ms of CPU time, and 50ms of real time, in part because the double-and-add method sketched in proposition~\ref{pro:Kamal_log} is particularly efficient in~$2$-power level.

We then generate an~$\F_\l$-basis of the subspace of~$J_1(16)(\F_{23})[5]$ that affords~$\rho_{f,\l}$ by using strategy~\ref{algo:Strategy_CycloExp}. This takes 710ms of CPU time, and 220ms of real time. On our way, we confirm that the rational canonical form of~$\rho_{f,\l}(\Frob_{23})$ is~$\smat{0}{-2}{1}{2} \in \GL_2(\F_\l)$; this was the only possibility, since during the first step, we had determined from the value of~$a_{23}(f) \bmod \l$ that the characteristic polynomial of~$\rho_{f,\l}(\Frob_{23})$ is~$x^2+2x+2$, which is separable mod~$5$.

We must now lift a basis of the representation space to~$J_1(16)(\Z_{23^4}/23^{7})[5]$. Actually, since we know now that the action of~$\Frob_{23}$ on the representation space is cyclic, we can afford to only lift one~$5$-torsion point, and then recover a basis by applying~$\Frob_{23}$ to it. This lifting takes 460ms of CPU time, and 260ms of real time.

Then, we generate all the points of the representation space over~$\Z_{23^4}/23^{7}$ by mixing the group law of the Jacobian and the action of~$\Frob_{23}$, and we evaluate the resulting points by 20 versions of the evaluation map~$\alpha$. All this takes 380ms of CPU time, and 80ms of real time.

Finally, we compute the corresponding 20 versions of the polynomial~$F(x)$, and keep the nicest one. This takes 7ms of CPU time, and 2ms of real time.

In the end, we find that our Galois representation is described by the polynomial
\begin{align*} F(x) = \ & x^{24} - 18x^{23} + 144x^{22} - 682x^{21} + 2141x^{20} - 4908x^{19} + 9014x^{18} \\ -& 14032x^{17} + 18606x^{16} - 20928x^{15} + 20086x^{14} - 15568x^{13} + 9009x^{12} \\  -& 5122x^{11} + 3206x^{10} - 1778x^9 + 5384x^8 - 9242x^7 + 7866x^6 - 4818x^5 \\ +& 1613x^4 - 124x^3 - 28x^2
+ 4x - 2 \in \Q[x]. \end{align*}
The whole computation took about 2.4s of CPU time, and 700ms of real time.

\begin{rk}
Since the computation also returns an indexation of the~$23$-adic roots of~$F(x)$ by the nonzero vectors of~$\F_\l^2$, we can easily compute a polynomial describing the projective version if we wish to do so, by gathering symmetrically (e.g. summing) the roots of~$F(x)$ along the vector lines of~$\F_\l^2$. We find the polynomial
\[ x^6 - 18x^5 + 120x^4 - 400x^3 + 680x^2 - 208x - 896 \in \Q[x], \]
which has one rational root (at $x=8$) and one irreducible factor of the degree~$5$. The representation~$\rho_{f,\l}$ is thus reducible, a fact that can easily be checked independently.
\end{rk}

\begin{rk}\label{rk:Qf_bad_1}
Experimenting shows that if we had used evaluation maps from the Jacobian to~$\A^1$ based on rational forms instead of~$q$-expansions (cf. section~\ref{sect:Eval}, Strategy 2), we would have had to increase the~$p$-adic accuracy to about~$O(23^{600})$, which would have slowed down the computation by a factor of about~30.
\end{rk}

\subsubsection*{$\Delta$ mod 13}

As a second example, we compute the representation attached to 
\[ \Delta = \mfref{1.12.a.a} = q - 24q^2 +252 q^3 +O(q^4) \]
mod~$\l = 13$. By the arguments presented in section~\ref{sect:rho_in_XH}, this representation~$\rho_{\Delta,13}$ is found in the~$13$-torsion of the Jacobian~$J_1(13)$ of the modular curve~$X_1(13)$, whose genus is again~$2$.

Since we know that the image of the representation is going to be the whole of~$\GL_2(\F_{13})$, this time we look for a prime~$p$ up to~$200$. This turns out not to be necessary: indeed, for $p=73$, the order of the image of the Frobenius is again~$a=4$ only. However, this whole search only took 720ms of CPU time, and 110ms of real time. Anyway, the computation proceeds with~$p=73$.

We choose to work at accuracy $O(73^{44})$, so as to be able to identify rationals of height up to $10^{40}$. Constructing a~$73$-adic Makdisi model of~$X_1(13)$ with residue degree~$a=4$ at this accuracy takes 190ms of CPU time, and 80ms of real time. This includes spotting the elliptic curve
\[ y^2 = x^3+25x+36, \]
which has all its~$13$-torsion defined over the degree-4 unramified extension of~$\Q_{73}$.

We then generate an~$\F_{13}$-basis of the subspace of~$J_1(13)(\F_{73})[13]$ that affords~$\rho_{\Delta,13}$. This takes 1030ms of CPU time, and 460ms of real time. On our way, we confirm that the rational canonical form of~$\rho_{\Delta,13}(\Frob_{73})$ is~$\smat{0}{-5}{1}{6} \in \GL_2(\F_{13})$, which as in the previous example we already knew from the first step.

Lifting a~$13$-torsion point which generates the representation space under~$\Frob_{73}$ to accuracy~$O(73^{44})$ takes 2.2s of CPU time, and 940ms of real time. After this, generating all the points of the representation space over~$\Z_{73^4}/73^{44}$ and evaluating them takes 6.8s of CPU time, and 970ms of real time.

Finally, we compute 24 versions of the polynomial~$F(x)$ and keep the nicest one, which takes 360ms of CPU time, and 60ms of real time.

In the end, we find that our Galois representation is described by a polynomial of the form
\[ x^{168} + \frac{290398}{10103}x^{167} + \cdots - \frac{36719}{10103} \in \Q[x], \]
whose coefficients have~$10103^2$ as a common denominator, and numerators of up to nearly 40 decimal digits.

The whole computation took about 11.3s of CPU time, and 2.6s of real time. As a comparison, a few years ago, the computation~\cite{algo} of the same representation by the complex-analytic method took about 5 minutes of \emph{real time} on the supercomputing cluster~\cite{Plafrim}, even though we parallelised it over dozens of cores.

\begin{rk}\label{rk:Qf_bad_2}
Experimenting shows that if we had used evaluation maps from the Jacobian to~$\A^1$ based on rational forms instead of~$q$-expansions (cf. section~\ref{sect:Eval}), we would have had to increase the~$p$-adic accuracy to about~$O(73^{8000})$.
\end{rk}

\subsubsection*{$\Delta$ mod 19}

We now try a larger example, that of the representation attached to~$\Delta$ mod~$\l=19$, which is found in the 19-torsion of the Jacobian of a curve of genus~$g=7$, namely~$X_1(19)$.

After having tried all primes~$p \leqslant 1000$, we select~$p=107$, since it allows us to work in residue degree~$a=6$. The search took 6s of CPU time, and 1s of real time.

We choose to work at accuracy $O(19^{247})$, so as to be able to identify rationals of height up to $10^{250}$. Constructing a~$107$-adic Makdisi model of~$X_1(19)$ with degree 6 at this accuracy takes in 11s of CPU time and 4.7s of real time.

After this, generating a basis of the representation space over~$\F_{107}$ takes 39s of CPU time and 12s of real time.

Next, lifting a 19-torsion point to accuracy~$O(107^{247})$ took 19 minutes of CPU time and 2m54s of real time, after which generating and evaluating all the other points took 32m30s of CPU time and 4m30s of real time. Finally, the computation of 12 versions of~$F(x)$ took 5s of CPU time and under 1s of real time.

In total, the computation took under 1h of CPU time, and under 8m of real time. In comparison, a few years ago, the computation~\cite{algo} of the same representation by the complex-analytic method took about 40 minutes of \emph{real time} on the supercomputing cluster~\cite{Plafrim}, even though we parallelised it over dozens of cores. This difference, although still impressive, is less striking than in the previous example, because we have to compute~$\ell^2-1=360$ torsion points whereas the author's laptop only has 4 cores, but also because we had to work in slightly higher residual degree this time.

\subsubsection{``Larger'' examples}

We now demonstrate the performance of our method on ``larger'' examples, which we run on the supercomputing cluster~\cite{Plafrim}. Since we perform parallel computations there using the MPI threading engine, we are no longer able to accurately measure the CPU times, and only give real times from now on.

\subsubsection*{\mfref{7.8.a.a} mod 13}

The following example was executed on 64 cores. Let 
\[ f = \mfref{7.8.a.a} = q - 6q^{2} - 42q^{3}+ O(q^{4}) \]
be the unique newform of weight~$k=8$ and level~$\Gamma_0(7)$ having rational coefficients, and let~$\ell=13$. The representation~$\rho_{f,13}$ is found in the~$13$-torsion of the Jacobian of the modular curve~$X_1(7 \cdot 13)$. This curve has genus~$265$, which is far too high for our method to apply; however, the arguments presented in section~\ref{sect:rho_in_XH} show that~$\rho_{f,13}$ actually occurs in the Jacobian of a curve~$X_H(7 \cdot 13)$ of genus~$g=13$ only. Therefore, our implementation chooses to use this curve to compute this representation.

We tried all the primes~$p \leqslant 1000$, and selected~$p=239$ since it lets us work in residue degree~$a=4$. The search took 4s.

Next, we generated a~$239$-adic Makdisi model of~$X_H(7 \cdot 13)$ with accuracy~$O(239^{256})$ in 37s.

After this, we computed a basis of the representation space over~$\F_{239}$ in 1m35s, lifted one of its points to accuracy~$O(239^{256})$ in 6m30s, computed and evaluated all the points in the representation space in 2m10s, and generated and selected a version of the polynomial~$F(x)$ in 200ms.

In total, the computation took 11m15s on 64 cores. In comparison, a few years ago, the computation~\cite{companion} of the same representation by the complex-analytic method took a little more than half a day (also of real time) on the Warwick mathematics institute computing cluster, also on 64 cores.

\subsubsection*{\mfref{5.6.a.a} mod 13}\label{ex:Tp}

Let 
\[ f = \mfref{5.6.a.a} = q +2q^{2} - 4q^{3}+ O(q^{4}) \]
be the unique newform of weight~$k=6$ and level~$\Gamma_0(5)$. The representation $\rho_{f,13}$ occurs with multiplicity $1$ in the $13$-torsion of the Jacobian of the genus~$13$ modular curve $X_H(5 \cdot 13)$; however, we observed in example~\ref{ex:Frobp_cannot_cut} that the $\F_{13}$-subspace $T_{f,13}$ of $J_H(5 \cdot 13)[13]$ which affords $\rho_{f,13}$ cannot be isolated by the action of $\Frob_p$ for any prime $p$. We therefore use this example to illustrate the variant of our method presented in section~\ref{sect:Tp}, again on 64 cores on~\cite{Plafrim}.

Indeed, a search over the primes $p \leqslant 1000$ with the methods presented in that section reveals that $T_{f,13}$ may be isolated as
\[ T_{f,13} = \ker\big(T_p - a_p(f) \vert_{ J_H(5 \cdot 13)[13]} \big) \]
for many primes $p$. This is in particular the case for $p=439$, which has the extra advantage that the order of $\rho_{f,13}(\Frob_p)$ is $a=6$ only. This search takes 3.7s.

As we wish to be able to identify rationals of height up to $10^{300}$, our implementation proceeds by constructing a $439$-adic Makdisi model of $J_H(5 \cdot 13)$ over $\Q_{439^6}$ with accuracy $O(439^{228})$, which takes 40s.

Our implementation then generates a few random points of $J_H(5 \cdot 13)[13](\F_{439^6})$ by the method outlined in section~\ref{sect:cycloexp}. The first of these points, which was generated in the subspace killed by $\Phi_3(\Frob_{439})$, spans an $\F_{13}[T_{439}]$-module of $\F_{13}$-dimension $2$, on which the matrix of $T_{439}$ is
\[ \left[ \begin{matrix} 0 & 0 \\ 1 & 8 \end{matrix} \right]. \]
Since $a_{439}(f) = 8 \bmod 13$, this module does not yet contain $T_{f,13}$, so we enlarge it by including the second random $13$-torsion point, which was generated in the subspace killed by $\Phi_2(\Frob_{439})$. The dimension of the $\F_{13}[T_{439}]$-module spanned by these two points is now 3, and the matrix of $T_{439}$ is now
\[ \left[ \begin{matrix} 0 & 0 & 0 \\ 1 & 8 & 0 \\ 0 & 0 & 8 \end{matrix} \right], \]
so this time we can extract an $\F_{13}$-basis of $T_{f,13}$ from this module. All this takes~4m10s.

From this point on, we proceed as usual. Lifting a generator of the $\F_13[\Frob_{439}]$ to accuracy~$O(439^{228})$ takes 15 minutes, computing and evaluating all the points in $T_{f,13}$ takes 3m40s, and generating $2$ versions of the polynomial~$F(x)$ and selecting the nicest one takes 180ms.

In total, the computation took 24 minutes on 64 cores.

\subsubsection*{$\Delta$ mod 29}

As a last example, we compute the representation attached to~$\Delta$ mod~$\l=29$. The smallest curve (that we know of) whose Jacobian contains this representation is~$X_1(29)$, whose genus is~$g=22$. Our implementation thus used this curve to compute~$\rho_{\Delta,29}$, again on~\cite{Plafrim} but using two machines with 42 cores each this time.

We tried all primes~$p \leqslant 1000$, and decided to use~$p=191$ because it lets us work in residual degree~$a=4$ only. This search took 1.3s.

Next, we generated a~$191$-adic Makdisi model of~$X_1(29)$ with accuracy~$O(191^{2048})$ in 21m.

After this, we computed a basis of the representation space over~$\F_{191}$ in 12m, lifted one of its points to accuracy~$O(191^{2048})$ in 6h10m, computed and evaluated all the points in the representation space in 4h10m, and generated and selected a version of the polynomial~$F(x)$ in 2m.

In total, the computation took a little less than 11h. In comparison, a few years ago~\cite{algo}, the computation (also on~\cite{Plafrim}) of the same representation by the complex-analytic method took about 3 days, even though it used about twice as many cores.

\begin{rk}
These examples show that the determination of an optimal~$p$-adic Makisi model of the modular curve is very far from being the bottleneck of the computation of a Galois representation. Besides, the last example also demonstrates that our~$p$-adic lifting method~\cite{Hensel} remains reasonably efficient in high genera.
\end{rk}

\subsection{Comparison with the complex-analytic method}\label{sect:compare}

The previous examples show that our implementation of the~$p$-adic method significantly outperforms our implementation of the complex-analytic method. That we wrote the former in C language whereas the latter was written in Python probably plays a part in this, but there are other, more fundamental reasons for this difference of performance.

Indeed, in order to generate~$\ell$-torsion points, the complex-analytic method begins by computing a high-accuracy approximation over~$\C$ of a period lattice of the modular curve, which takes a significant amount of time since it requires in particular computing the~$q$-expansion of a basis of the space of cusp forms of weight~2 to high accuracy. On the contrary, the~$p$-adic lifting method starts in ``low accuracy'', that is to say mod~$p$, where torsion points can be obtained easily thanks to fast exponentiation; therefore it does not suffer from this overhead. This explains in particular the major performance differences observed with the ``small'' examples above; thanks to the~$p$-adic approach, these calculations can now be executed on a personal laptop in very reasonable time.

Besides, as explained in~\cite[6.4]{Hensel}, since the evaluation map from the Jacobian to~$\A^1$ is by design Galois-equivariant, the~$p$-adic method can save a lot of effort by computing and evaluating not all the points of the representation space, but only one per orbit under the Frobenius~$\Frob_p$. In contrast, the complex method can only use complex conjugation, which has order~2 and thus can only halve the amount of work.

Another pleasant feature of the~$p$-adic approach is that it can naturally deal with non-squarefree levels, as demonstrated by the first example above which took place in level~$N=16$. On the contrary, non-squarefree levels are problematic for the complex method, since the computation of the periods of the modular curve requires the determination of Atkin-Lehner pseudo-eigenvalues~\cite[3.2.3]{companion}, which cannot in general be easily read off the coefficients of a newform when the level is not squarefree~\cite[2.1.2]{companion}.

\bigskip

Our implementation is available on the GitHub repository~\cite{Github} and in a development branch of~\cite{gp}.

\bigskip

\section*{Acknowledgements}
\addcontentsline{toc}{section}{Acknowledgements}

The author wishes to express him warm gratitude to Kamal Khuri-Makdisi for his algorithm to compute in Jacobians and his construction of moduli-friendly modular forms, which both play a central role in the methods discussed in this article, and for providing a detailed explanation of the algorithm outlined in proposition~\ref{pro:Kamal_log}; and to Bill Allombert for his help with parallel computation on the computing cluster~\cite{Plafrim} using the~\cite{gp} C language library, without which the large-scale computations shown in section~\ref{sect:Examples} would not have been possible.

Experiments presented in this paper were carried out using the PlaFRIM experimental testbed, supported by Inria, CNRS (LABRI and IMB), Universit\'e de Bordeaux, Bordeaux INP, and Conseil R\'egional d’Aquitaine (see \url{https://www.plafrim.fr/}).

\end{document}